\documentclass[11pt]{article}

\usepackage{amsmath,amssymb,amsfonts,amsthm}
\usepackage{moreverb,rotating,graphics}
\usepackage[T1]{fontenc}
\usepackage[latin1]{inputenc}
\usepackage{epsfig}
\usepackage[francais,english]{babel} 
\usepackage{color}
\usepackage{stmaryrd}
\usepackage{dsfont}
\usepackage{float}
\usepackage{setspace}

\newtheorem{theorem}{Theorem}
\newtheorem{Prop}[theorem]{Proposition}
\newtheorem{remark}{Remark}

\newtheorem{lemma}[theorem]{Lemma}

\def\NN{\mathbb{N}}
\def\ZZ{\mathbb{Z}}
\def\RR{\mathbb{R}}

\numberwithin{equation}{section}

\title{A numerical approach related to defect-type theories for some weakly random problems in homogenization} 
\author{A. Anantharaman$^{1}$ \thanks{The author acknowledges EADS IW for financial support.} \, and C. Le Bris$^{2}$\\ \\
       \footnotesize{Universit\'e Paris-Est, CERMICS, Project-team
         Micmac, INRIA-Ecole des Ponts,} \\  
\footnotesize{6 \& 8 avenue Blaise Pascal,
         77455 Marne-la-Vall\'ee Cedex 2, France} \\ 
\footnotesize{$^{1}$ \tt ananthaa@cermics.enpc.fr, \; $^{2}$ lebris@cermics.enpc.fr}\\
}

\begin{document}

\maketitle


\begin{abstract}
\noindent
We present in this paper an approach for computing the  homogenized
behavior of a medium that is a small random perturbation of a periodic
reference material. The random perturbation we consider is, in a sense
made precise in our work,  a rare event at the microscopic level. It
however affects  the macroscopic properties of the material, and we
indeed provide a method to compute the first and second-order
corrections. To this end, we formally establish
an asymptotic expansion of the macroscopic properties. Our perturbative
approach shares
common features with a defect-type theory of solid state physics. The
computational efficiency of the approach is demonstrated.\\
\end{abstract}

{\footnotesize{{\bf Keywords}: Homogenization; Random Media; Defects}\\

}

\section{Introduction}

Composite materials are increasingly used in  industry. For instance,
modern aircrafts consist, for more than 50$\%$,  of composite
materials. Generally speaking, composites are heterogeneous materials
obtained by mixing two phases, a matrix and reinforcements (or
inclusions). When appropriately designed, these materials outperform traditional
 materials, notably because they combine robustness and lightness. Their
 use however raises new challenges. The behavior of these materials
 under extreme conditions has to be predicted carefully, so as to avoid,
 in the worst case scenario,  separation of the components (think of a plane hit by
 thunder). While it is possible to create an infinity of composites
 starting from the same elementary components,  it is out of question to
 actually construct and experimentally test each and every possible
 combination. Characterizing {\it a priori} the properties of a given
 composite material, not yet synthetized or assembled, is therefore instrumental.\\

A brute force numerical approach, consisting in directly solving the
classical boundary value
problems modelling the behavior of the material, is not practical. The
heterogeneities indeed often occur at a  scale  $\epsilon$ much finer than the
overall typical lengthscale (say,~1) of the material itself. A finite
element  mesh would, for example, need to be of size less than
$\epsilon$ in order  to capture the correct behavior.  The number of degrees of
freedom 
would then be proportional to $\epsilon^{-d}$ (where $d$ denotes the
dimension of the ambient physical space) and would yield, for $\epsilon$
small, a  heavy computational cost one cannot necessarily afford.\\       

The aim of homogenization is to provide a practical alternative to the brute
force numerical approach. In a nutshell, homogenization consists in
replacing a possibly complicated heterogeneous material with a
homogeneous material sharing  the same macroscopic properties. 
It allows for eliminating  the fine scale, up to an error which is
controlled by $\epsilon$, the size of this fine scale as compared to the
macroscopic size. Homogenization is a well-established theory (see \cite
{JKO} for a comprehensive textbook), which, in a simplified picture,  can
be seen as averaging partial differential equations that have highly-oscillating coefficients.\\ 

Of course, the structure of the material, and more precisely the way the
constituents are combined, have a deep influence on the results of the
homogenization process. The simplest possible situation is the periodic
situation. At the fine scale, a unit cell is repeated in a periodic manner in all
directions. Then, in simple cases (say, to be schematic and to fix the ideas, linear
well-posed equations), the homogenized material is characterized only
using the solution of simple problems on the unit cell, called the cell
problems. The role of these cell problems is to encode the information
of the micro-scale and
convey it to the macro-scale. Related cases, such as pseudo-periodic
materials, can be treated similarly.\\

\begin{figure}[h]\label{composite}
\center
\includegraphics[width = 10cm, height = 7cm]{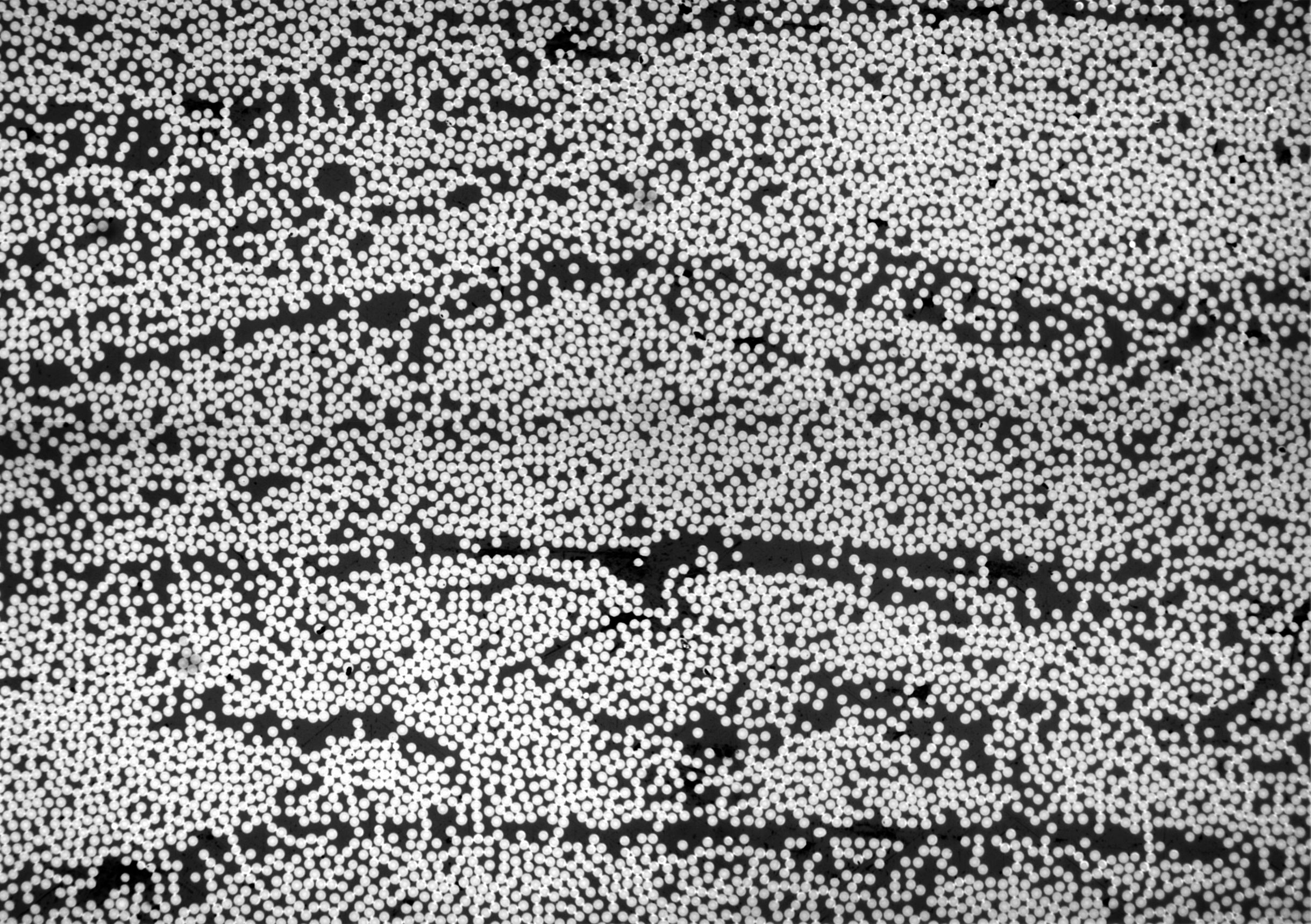}
\caption{Two-dimensional cut of a composite material used in the aeronautics industry, extracted from \cite{Matthieu} and reproduced with permission of the author. It is clear that this material is not periodic, yet there is some kind of an underlying periodic arrangement of the fibers.}
\end{figure}

As Figure \ref{composite} shows, real life materials are however not
often periodic. In particular because of uncertainties and flaws in the
industrial process,  composites often do not exhibit a perfect periodic structure,
even though it was the original plan. A suitable way to account for this
is to use random modelling.  Although the mathematical theory for
homogenization of random materials under classical assumptions
(ergodicity and stationarity) is well known, the practice is quite
involved. The cell problems are defined over the whole space
$\mathbb{R}^d$ and not simply on a ``unit'' cell. The numerical
approximation of such problems using Monte-Carlo type computations is incredibly costly: the cell problems are
truncated on a bounded domain, many possible
realizations of the materials are considered,  averages are
performed. Consequently, in the context of random modelling, the
benefits of homogenization over the direct attack of the original
composite material are arguable.\\

Our line of thoughts, and the approach we try to advocate here, are based
on the following two-fold observation: classical random homogenization is
costly but perhaps, in a number of situations, not necessary. A more
careful examination of  Figure \ref{composite} indeed shows that albeit
not periodic, the material is not totally random. It may probably be fairly considered as a perturbation of a periodic material. The homogenized
behavior should expectedly be close to  that of the underlying periodic
material, up to a small error  depending on the amount of randomness present.\\

The aim of this paper is to give a practical example of theory following
the above philosophy. We introduce and study a specific model for such a randomly perturbed periodic material, which we also call a \emph{weakly random material}. More precisely, we are interested in the homogenization of the following elliptic problem
\begin{eqnarray*}
\left \{
\begin{aligned}
&-\mathrm{div} \left((A_{per}(\frac{x}{\epsilon}) + b_{\eta}(\frac{x}{\epsilon}, \omega) C_{per}(\frac{x}{\epsilon}) ) \nabla u_{\epsilon}\right) = f(x) \; \mathrm{in} \; \mathcal{D} \subset{\mathbb{R}^d}, \\
&u_{\epsilon} = 0 \; \mathrm{on} \; \partial \mathcal{D}.
\end{aligned}
\right.
\end{eqnarray*}
Here the tensor $A_{per}$ models a reference $\mathbb{Z}^d$-periodic
material which is randomly perturbed by the $\mathbb{Z}^d$-periodic
tensor $C_{per}$, the stochastic perturbation being encoded in the
stationary ergodic scalar field $b_{\eta}$. In the present work, the law
of the random variable $b_{\eta}(x,\cdot)$ is a Bernoulli distribution with
parameter $\eta$ (that is, $b_{\eta}$ is equal to $1$ with probability
$\eta$ and $0$ with probability $1-\eta$. Using an asymptotic analysis in
terms of $\eta$, we will develop an homogenization theory
for  $A_{\eta}\left(\frac{x}{\epsilon}, \omega\right) = A_{per}\left(\frac{x}{\epsilon}\right)
+ b_{\eta}\left(\frac{x}{\epsilon}, \omega\right) C_{per}$ based on the
similar theory for  $A_{per}\left(\frac{x}{\epsilon}\right)$.\\

In short, let us say that the main result of this article is to formally derive an expansion 
$$A_{\eta}^* = A_{per}^* + \eta \bar{A_1^*} + \eta^2 \bar{A_2^*} + o(\eta^2),$$
where $A_{\eta}^*$ and $A_{per}^*$ are the homogenized tensors associated with $A_{\eta}$ and $A_{per}$ respectively. The first-order correction $\bar{A_1^*}$ is obtained as the limit, when $N \rightarrow \infty$, of a sequence of tensors $A_{1}^{*,N}$ computed on the supercell $[-\frac{N}{2},\frac{N}{2}]^d$. It is the purpose of Proposition \ref{conva1} to prove the convergence of $A_{1}^{*,N}$. The second-order term $\bar{A_2^*}$ is likewise defined as a limit, this time up to extraction, of a sequence of tensors $A_{2}^{*,N}$ when $N \rightarrow \infty$. The proof of the boundedness of the sequence $A_{2}^{*,N}$ which implies this convergence up to extraction is not given here for it relies on long and technical calculations. Actually, we strongly believe that $A_{2}^{*,N}$ is a convergent sequence and that we can write the expression of the limit. We refer the reader to \cite{these_Arnaud} for the details. We also stress that these corrections are achieved through purely deterministic computations.\\ 

 The above
setting  is of course \emph{one}
possible setting where we may develop our theory, but not the only one. More general distributions are
studied in \cite{ALB2}. Other forms of random perturbations of periodic
problems, in the spirit of \cite{BLL}, could also be addressed. Moreover, we
have deliberately considered the simplest possible equation (a scalar,
linear second order elliptic equation in divergence form) to avoid any
unnecessary technicalities and fundamental difficulties. Other equations could be
considered, although it is not currently clear (to us, at least) how
general  our theory is in this respect.
\\

With the ideas developped here (and originally introduced and further
mentioned in
\cite{ALB2,ALBcras,enumath}), we work in the footsteps of
many previous contributors who have considered perturbative approaches
in homogenization. In \cite{Tartar} and \cite{AG}, a deterministic setting in which an asymptotic expansion is assumed on the properties of the material (the latter being not necessarily periodic) is studied under the name ``small amplitude homogenization''. In
\cite{Sakata}, the case of a Gaussian perturbation with a small variance is addressed from a
mechanical point of view. Our setting here is particular, because  our random
perturbation has order one in amplitude. It is only \emph{in law} that the
perturbation  considered is small. The corrections obtained are
therefore intrinsically   different from those obtained in other
settings (including settings we ourselves consider elsewhere, see
\cite{these_Arnaud, ALB2}). Also, the
present perturbative theory has unanticipated close connections
with some classical defect-type theories used in solid state physics.\\

We emphasize  that, contrary to what is presented in a companion paper~\cite{ALB2} for some 
other distributions, the theoretical results we obtain below in the Bernoulli case
are only formal. We are unfortunately unable to fully justify our
manipulations except in the one-dimensional case. Nevertheless, we can
prove that the terms we obtain as first-order and second-order
corrections are indeed finite and well defined. Our numerical results, on the other hand, show the efficiency of the approach. They
somehow constitute a proof of the definite validity of our perturbative approach, although we wish to remain
cautious. Note that due to the prohibitive cost of three-dimensional random homogenization problems and the limited computing facilities we have access to, our tests are performed in dimension two.\\ 

This paper is organized as follows. For the sake of consistency and the
reader's convenience, we start by recalling in Section 2 some classical
results of periodic and stochastic elliptic homogenization. Then we
introduce our perturbative model in Section 3, and explain how we obtain
the first-order and second-order correction by means of an ergodic
approximation. Our elements of proof are exposed in Section~3.
Our two-dimensional numerical tests are presented  in Section 4.  The
appendix contains explicit computations in the one-dimensional case as well as some useful technical lemmas.\\  

Throughout this paper, and unless otherwise mentioned, $K$ denotes a constant that depends at most on the ambient dimension $d$, and on the tensors $A_{per}$ and $C_{per}$. The indices $i$ and $j$ denote indices in $\llbracket 1, d \rrbracket$.


\section{Some classical results of elliptic homogenization} \label{classic}

We recall here some classical well-known results regarding linear elliptic periodic and stochastic homogenization. The reader familiar with homogenization theory can easily skip this section and directly proceed to Section 3.

\subsection{Periodic homogenization}

Consider $A$ a $\ZZ^d$-periodic tensor field from $\RR^d$ to $\RR^{d\times d}$ , that is
$$\forall k \in \ZZ^d, A(x+k) = A(x) \; \; \mathrm{almost \; everywhere \; in \;} x \in \RR^d. $$
We assume that $A \in L^{\infty}(\RR^d, \RR^{d\times d})$ and $A$ is coercive, which means that there exist $\lambda >0 $ and $\Lambda >0 $ such that
\begin{eqnarray} \label{boundaper}
\forall \xi \in \RR^d, \mathrm{\; a.e \; in \;} x \in \RR^d, \; \lambda |\xi|^2 \leq A(x)\xi\cdot \xi \; \mathrm{and} \; |A(x) \xi| \leq \Lambda |\xi|. 
\end{eqnarray}
Consider now a material occupying a bounded domain $\mathcal{O} \subset \RR^d$. The constitutive properties of this material are supposed to be periodic, the scale of periodicity being $\epsilon$, and we assume that these properties are given by the tensor $A_{\epsilon}(x) = A \left( \frac{x}{\epsilon} \right)$.\\

We consider the following canonical elliptic problem: $f \in L^2(\mathcal{O})$ being given, find $u_{\epsilon} \in H^1_0(\mathcal{O})$ solution to
\begin{eqnarray} \label{homper}
\left \{
\begin{aligned}
&-\mathrm{div} \left(A_{\epsilon} \nabla u_{\epsilon}\right) = f \; \mathrm{in} \; \mathcal{O}, \\
&u_{\epsilon} = 0 \; \mathrm{on} \; \partial \mathcal{O}.
\end{aligned}
\right.
\end{eqnarray}    

A direct numerical handling of (\ref{homper}) using finite elements has a heavy computational cost since the scale $\epsilon$ of the heterogeneities requires a fine mesh. The aim of homogenization is to take the limit $\epsilon \rightarrow 0$ in (\ref{homper}) so as to replace the heterogeneous material with a homogeneous material. To this end, let us define the periodic cell problems on the unit cell $Q = [-\frac{1}{2},\frac{1}{2}]^d$ by: $\forall i \in \llbracket 1, d \rrbracket,$ 
\begin{eqnarray} \label{cellper}
\left  \{
\begin{aligned}
& -\mathrm{div}\left(A(\nabla w_{i} + e_i ) \right) = 0 \quad \mathrm{in} \, \, Q, \\
& w_{i} \; \; \mathbb{Z}^d-\mathrm{periodic},  
\end{aligned}
\right.
\end{eqnarray}
where $e_i$ is the $i$-th canonical vector of $\RR^d$. Problem (\ref{cellper}) has a solution unique up to the addition of a constant. Note that the number of cell problems is equal to the dimension of the space.\\

The homogenized tensor $A^*$ is then given by: 
\begin{eqnarray} \label{homtens}
\forall (i,j) \in \llbracket 1,d \rrbracket^2, \; \;A^*_{ji} = \int_Q A (\nabla w_i + e_i) \cdot e_j.
\end{eqnarray}
Using (\ref{cellper}), it also holds
$$A^*_{ji} = \int_Q A (\nabla w_i + e_i) \cdot (\nabla w_j + e_j).$$
Notice that in this periodic setting $A^*$ is a constant matrix.\\

Finally, let us define the homogenized solution $u_0$ as the unique solution in $H^1_0(\mathcal{O})$ to 
\begin{eqnarray} \label{homsol}
\left \{
\begin{aligned}
&-\mathrm{div} \left(A^* \nabla u_{0}\right) = f \; \mathrm{in} \; \mathcal{O}, \\
&u_{0} = 0 \; \mathrm{on} \; \partial \mathcal{O}.
\end{aligned}
\right.
\end{eqnarray}
Solving (\ref{cellper}) and (\ref{homsol}) is much simpler than directly solving (\ref{homper}) for the fine scale $\epsilon$ has disappeared. It is well-known (see \cite{JKO} for instance) that 
\begin{eqnarray} \label{convl2}
u_{\epsilon} \underset{\epsilon \rightarrow 0}{\rightarrow} u_0 \quad \mathrm{in} \; L^2(\mathcal{O})
\end{eqnarray}
and 
\begin{eqnarray} \label{convh1}
u_{\epsilon} - u_0 - \epsilon \sum_{i=1}^{d} w_i\left(\frac{\cdot}{\epsilon}\right)\frac{\partial u_0}{\partial x_i}  \underset{\epsilon \rightarrow 0}{\rightarrow} 0\quad \mathrm{in} \; H^1(\mathcal{O}).
\end{eqnarray}
The functions $w_i$ are also called the correctors, since they allow for the \emph{strong} convergence in (\ref{convh1}). Convergences (\ref{convl2}) and (\ref{convh1}) show the relevance of the homogenization process: $u_{\epsilon}$ can be replaced by $u_{0}$ or more accurately $u_{0} + \epsilon \sum_{i=1}^{d} w_i\left(\frac{x}{\epsilon}\right)\frac{\partial u_0}{\partial x_i}(x)$, which are easier to compute.

\subsection{Stochastic homogenization}

Throughout the article, $(\Omega , \mathcal{F}, \mathbb{P})$ denotes a probability space, $\mathbb{P}$ the probability measure and $\omega \in \Omega$ an event. We denote by $\mathbb{E}(X)$ the expectation of a random variable $X$.\\

We assume that the group $(\ZZ^d,+)$ acts on $\Omega$ and denote by $\tau_k, k \in \mathbb{Z}^d,$ the group action. We also assume that this action is measure-preserving, that is,
$$\forall \mathcal{A} \in \mathcal{F},  \forall k \in \mathbb{Z}^d, \; \mathbb{P}(\mathcal{A}) = \mathbb{P}(\tau_k \mathcal{A}),$$
and ergodic: 
$$\forall \mathcal{A} \in \mathcal{F},  (\forall k \in \mathbb{Z}^d, \mathcal{A} = \tau_k \mathcal{A})  \implies (\mathbb{P}(\mathcal{A})=0 \; \mathrm{or} \; \mathbb{P}(\mathcal{A})=1).$$

We call $F \in L^1_{loc}(\RR^d, L^1(\Omega))$ stationary if 
\begin{eqnarray}\label{statio}
 \forall k \in \ZZ^d, \; F(x+k, \omega) = F(x, \tau_k \omega)\; \; \mathrm{almost \; everywhere \; in \;} x \in \RR^d \; \mathrm{and \; } \omega \in \Omega.
\end{eqnarray}

Notice that the notion of stationarity we use here is \emph{discrete}: the shifts in (\ref{statio}) are assumed to be integers. This is related to our wish to connect the random problems considered with some underlying periodic problems.
Notice also that for a deterministic $F$, stationarity amounts to $\ZZ^d$-periodicity.\\

Consider a stationary tensor field $A(x,\omega) \in L^{\infty}(\RR^d\times \Omega , \RR^{d\times d})$, such that (\ref{boundaper}) is almost surely satisfied by $A(\cdot,\omega)$, and a material occupying a bounded domain $\mathcal{O} \subset \RR^d$ modeled by $A\left(\frac{x}{\epsilon},\omega\right)$. 

We are interested in solving, for a deterministic function $f$,
\begin{eqnarray} \label{homstoch}
\left \{
\begin{aligned}
&-\mathrm{div} \left(A \left(\frac{x}{\epsilon}, \omega \right) \nabla u_{\epsilon}\right) = f(x) \; \mathrm{in} \; \mathcal{O}, \\
&u_{\epsilon} = 0 \; \mathrm{on} \; \partial \mathcal{O}.
\end{aligned}
\right.
\end{eqnarray}

In order to describe the behavior of $u_{\epsilon}$, we again need to define cell problems. Here they read (see \cite{JKO}):
\begin{eqnarray} \label{cellsol}
\left  \{
\begin{aligned}
& -\mathrm{div}\left(A(x,\omega)(\nabla w_{i} + e_i ) \right) = 0  \quad \mathrm{in} \, \, \RR^d,\\
& \nabla w_{i} \; \; \mathrm{stationary}, \quad \mathbb{E}\left(\int_{Q} \nabla w_{i}\right) = 0. \\ 
\end{aligned}
\right.
\end{eqnarray}

Problem (\ref{cellsol}) has a solution unique up to the addition of a (possibly random) constant.\\

Then we define the homogenized tensor $A^*$ by
\begin{eqnarray} \label{stochtens}
\forall (i,j) \in \llbracket 1,d \rrbracket^2, \displaystyle A^*_{ji} = \mathbb{E}\left(\int_Q A(y,\omega)(e_i  + \nabla_y w_i(y,\omega))\cdot e_j dy \right).
\end{eqnarray}
Notice that $A^*$ is deterministic and constant throughout the domain $\mathcal{O}$. The homogenized field $u_0$, which gives the asymptotic behavior of $u_{\epsilon}$ (in a sense similar to (\ref{convl2}) and (\ref{convh1})), is also deterministic. It is the unique solution in $H^1_0(\mathcal{O})$ to
\begin{eqnarray*}
\left \{
\begin{aligned}
&-\mathrm{div} \left(A^* \nabla u_{0}\right) = f \; \mathrm{in} \; \mathcal{O}, \\
&u_{0} = 0 \; \mathrm{on} \; \partial \mathcal{O}.
\end{aligned}
\right.
\end{eqnarray*}

The computation of the stochastic cell problems (\ref{cellsol}) is not an easy task since the problems are posed in an infinite domain ($\RR^d$) with a stationarity condition. As we have seen in the previous paragraph, when the material is periodic, the cell problems (\ref{cellsol}) reduce to the deterministic cell problems (\ref{cellper}) which are $\ZZ^d$-periodic and can thus be computed on the unit cell $Q$. Consequently, when the material under consideration is a stochastic perturbation of a reference periodic material, we expect the computation of the homogenized tensor to be tractable, up to an approximation. This is our motivation for proposing a perturbative approach.

\section{Homogenization of a randomly perturbed periodic material} \label{model}

\subsection{Presentation of the model}

In the stochastic framework (\ref{homstoch})-(\ref{cellsol})-(\ref{stochtens}), we now specifically consider the following tensor field in $\RR^d \times \Omega$:
\begin{eqnarray} \label{homogeta}
A_{\eta}(x,\omega) = A_{per}(x) + b_{\eta}(x,\omega)C_{per}(x).
\end{eqnarray}

Here $A_{per}$ and $C_{per}$ are two deterministic $\mathbb{Z}^d$-periodic tensor fields. Intuitively, $A_{per}$ is the reference material perturbed by $C_{per}$. The random character of the perturbation is encoded in the stationary ergodic scalar field $b_{\eta}$ , upon which we assume the expression
$$b_{\eta}(x,\omega) = \sum_{k \in \ZZ^d}\mathds{1}_{Q + k}(x) B_{\eta}^k(\omega),$$
where the $B_{\eta}^k$ are independent random variables having Bernoulli distribution with parameter $\eta$, meaning $B_{\eta}^k = 0$ with probability $1-\eta$ and $B_{\eta}^k = 1$ with probability $\eta$.

It is clear that as $\eta \rightarrow 0$ the perturbation becomes a rare event. However, the realization of this event modifies the microscopic structure of the material since it replaces, in a given cell, $A_{per}$ with $A_{per} + C_{per}$. \\

We additionally assume that there exist $0 < \alpha \leq \beta$ such that for all $\xi \in \RR^d$ and almost all $x \in \RR^d$,
\begin{eqnarray} \label{alphabound}
\alpha |\xi|^2 \leq A_{per}(x)\xi\cdot \xi, \qquad \alpha |\xi|^2 \leq \left(A_{per}+C_{per}\right)(x)\xi\cdot \xi,
\end{eqnarray}
\begin{eqnarray} \label{betabound}
|A_{per}(x) \xi| \leq \beta |\xi|, \qquad |\left(A_{per}+ C_{per}\right)(x)\xi| \leq \beta |\xi|.
\end{eqnarray}

We can therefore use for every $0 \leq \eta \leq 1$ the stochastic homogenization results recalled in Section \ref{classic}. The cell problems associated with (\ref{homogeta}) read, for $1\leq i \leq d$,
\begin{eqnarray} \label{etacell} 
\left  \{
\begin{aligned}
& -\mathrm{div}\left(A_{\eta}(\nabla w_{i}^{\eta} +e_i ) \right) = 0 \quad \mathrm{in} \, \, \RR^d, \\
& \nabla w_{i}^{\eta} \; \; \mathrm{stationary}, \; \; \mathbb{E}\left(\int_Q \nabla w_{i}^{\eta}\right) = 0,
\end{aligned}
\right.
\end{eqnarray}
and the homogenized tensor $A_{\eta}^{*}$ is given by
\begin{eqnarray} \label{homogal}
A_{\eta}^* e_i= \mathbb{E}\left(\int_Q A_{\eta}(\nabla w_{i}^{\eta} + e_i) \right), \; \; \mathrm{for \;}   1 \leq i \leq d.
\end{eqnarray}

Throughout the rest of this paper we denote by $w_i^0$ the solution to the $i$-th cell problem (\ref{cellper}) associated with $A_{per}$.\\

Because of the specific form of $A_{\eta}$, and more precisely because $A_{\eta}$ converges strongly to $A_{per}$ in $L^2(Q \times \Omega)$ as $\eta \rightarrow 0$, it is easy to see that:

\begin{lemma}
When $\eta \rightarrow 0$, $A_{\eta}^* \rightarrow A_{per}^{*}$. 
\end{lemma}

\begin{proof}
Fix $1 \leq i \leq d$. We start by proving that $\nabla w_i^{\eta}$ converges strongly in $L^2(Q \times \Omega)$ to $\nabla w_i^{0}$. Indeed, define $r_i^{\eta} = w_i^{\eta} - w_i^0$ solution to  
\begin{eqnarray} \label{celldiff}
\left  \{
\begin{aligned}
& -\mathrm{div}\left(A_{\eta}\nabla r_{i}^{\eta} \right) = \mathrm{div}\left(b_{\eta}C_{per}\left(\nabla w_i^0 + e_i \right) \right)\quad \mathrm{in} \, \, \RR^d,\\
& \nabla r_{i}^{\eta} \; \; \mathrm{stationary}, \quad \mathbb{E}\left(\int_{Q} \nabla r_{i}^{\eta}\right) = 0. \\ 
\end{aligned}
\right.
\end{eqnarray}

Standard cut-off and ergodicity arguments (see e.g the proof of Proposition 3.1 in \cite{BLL}) show that 
\begin{eqnarray*}
\|\nabla r_{i}^{\eta}\|_{L^2(Q \times \Omega)} &\leq& \frac{1}{\alpha} \|b_{\eta}C_{per}\left(\nabla w_i^0 + e_i \right)\|_{L^2(Q \times \Omega)}\\
&=& \frac{1}{\alpha} \|B_{\eta}^0\|_{L^2(\Omega)} \|C_{per}\left(\nabla w_i^0 + e_i \right)\|_{L^2(Q)}\\
&=& \frac{1}{\alpha} \sqrt{\eta} \|C_{per}\left(\nabla w_i^0 + e_i \right)\|_{L^2(Q)},
\end{eqnarray*}
where $\alpha$ is defined in (\ref{alphabound}), so that $\nabla w_i^{\eta} \underset{\eta \rightarrow 0}{\rightarrow} \nabla w_i^{0}$ in $L^2(Q \times \Omega)$ .\\

Next, it is straightforward to see that $A_{\eta}$ converges strongly to $A_{per}$ in $L^2(Q \times \Omega)$. We deduce from these two strong convergences that
$$A_{\eta}^* e_i = \mathbb{E}\left(\int_Q A_{\eta}(x,\omega) (\nabla w_{i}^{\eta} + e_i) \right) \underset{\eta \rightarrow 0}{\rightarrow} \int_Q A_{per}(\nabla w_{i}^{0} + e_i)=A_{per}^* e_i.$$
This concludes the proof.

\end{proof}

Our goal is now to find an asymptotic expansion for $A_{\eta}$ with respect to $\eta$ up to the second order.

\subsection{An ergodic approximation of the homogenized tensor}

We consider a specific realization $\tilde{\omega} \in \Omega$ of the tensor $A_{\eta}$ in the truncated domain $I_N = [-\frac{N}{2},\frac{N}{2}]^d$, with (for simplicity) $N$ an odd integer, and solve the following ``supercell'' problem:
\begin{eqnarray}\label{ergodiccorr} 
\left  \{
\begin{aligned}
& -\mathrm{div}\left(A_{\eta}(x,\tilde{\omega})(\nabla w_{i}^{\eta,N,\tilde{\omega}} + e_i ) \right) = 0 \quad \mathrm{in} \, \, I_N,\\
& w_{i}^{\eta,N,\tilde{\omega}} \, \, (N\ZZ)^d-\mathrm{periodic}. 
\end{aligned}
\right.
\end{eqnarray} 

Then an easy adaptation of Theorem 1 of \cite{BP}, stated in the continuous stationary setting, to our discrete stationary setting, shows that when $N$ goes to infinity,
\begin{eqnarray} \label{almost}
\displaystyle \frac{1}{N^d} \int_{I_N}A_{\eta}(x,\tilde{\omega})(\nabla w_{i}^{\eta,N,\tilde{\omega}}(x) + e_i ) dx \mathrm{\; converges \; to\;} A_{\eta}^*e_i \mathrm{\; almost \; surely \; in} \; \tilde{\omega} \in \Omega.
\end{eqnarray} Since $\displaystyle \frac{1}{N^d} \int_{I_N}A_{\eta}(x,\tilde{\omega})(\nabla w_{i}^{\eta,N,\tilde{\omega}}(x) + e_i ) dx$ is the tensor obtained by periodic homogenization of the tensor $A_{\eta}(x,\tilde{\omega})$ in the supercell $I_N$, it is also well-known (see \cite{JKO}) that the following bounds hold for all $(i,j) \in \llbracket 1,d \rrbracket^2$:
\begin{eqnarray*}
\frac{1}{N^d} \left(\int_{I_N}A_{\eta}^{-1}(x,\tilde{\omega})dx\right)^{-1} e_i \cdot e_j  \leq \frac{1}{N^d} \int_{I_N}A_{\eta}(x,\tilde{\omega})(\nabla w_{i}^{\eta,N,\tilde{\omega}}(x) + e_i )\cdot e_j dx \\
 \qquad \qquad \leq \frac{1}{N^d} \left(\int_{I_N}A_{\eta}(x,\tilde{\omega})dx \right) e_i \cdot e_j.
\end{eqnarray*}
 As a result, for all $N$ in $2\NN+1$, for all $0 \leq \eta \leq 1$ and almost all $\tilde{\omega}$ in $\Omega$, 
\begin{eqnarray}\label{Lebesgue}
\left|\frac{1}{N^d} \int_{I_N}A_{\eta}(x,\tilde{\omega})(\nabla w_{i}^{\eta,N,\tilde{\omega}}(x) + e_i )\cdot e_j dx \right| \leq \beta,
\end{eqnarray}
where $\beta$ is defined in (\ref{betabound}).
We then deduce from (\ref{almost}), (\ref{Lebesgue}) and the Lebegue dominated convergence theorem that
\begin{eqnarray} \label{expect}
\forall 1 \leq i \leq d, \; \; A_{\eta}^*e_i = \lim_{N \rightarrow + \infty} \frac{1}{N^d} \mathbb{E} \left(\int_{I_N}A_{\eta}(x,\omega)(\nabla w_{i}^{\eta,N,\omega}(x) + e_i ) \right)dx. 
\end{eqnarray}

\begin{remark}
A similar result holds for homogeneous Dirichlet and Neumann boundary conditions instead of periodic conditions in the definition (\ref{ergodiccorr}) of $w_{i}^{\eta,N,\tilde{\omega}}$ (see \cite{BP} for more details).
\end{remark}

Using now the fact that $b_{\eta}$ has a Bernoulli distribution in each cell of $\ZZ^d$, it is a simple matter to count the events and to make (\ref{expect}) more precise. We first define the set
\begin{eqnarray} \label{taun}
\mathcal{T}_N = \left\{ k \in \ZZ^d, Q+k \subset I_N \right\} = \left \llbracket -\frac{N-1}{2}, \frac{N-1}{2} \right \rrbracket^d.
\end{eqnarray}
The cardinal of $\mathcal{T}_N$ is of course $N^d$, and $\displaystyle \bigcup_{k \in \mathcal{T}_N} \{Q+k\} = I_N.$

We then have the following possible values for $A_{\eta}$: 
\begin{itemize}
\item $A_{\eta}(x,\tilde{\omega}) = A_{per}$ with probability $(1-\eta)^{N^d}$. 

In this case $w_{i}^{\eta,N,\tilde{\omega}} = w_{i}^{0}$ solves the usual periodic cell problem:
\begin{eqnarray*} 
\left  \{
\begin{aligned}
& -\mathrm{div}\left(A_{per}(\nabla w_i^0 + e_i)\right) = 0\quad \mathrm{in} \, \, Q, \\
& w_{i}^{0} \, \, \ZZ^d-\mathrm{periodic}. 
\end{aligned}
\right.
\end{eqnarray*} 

\item $A_{\eta}(x,\tilde{\omega}) = A_{per} + \mathds{1}_{\{Q+ k\}}C_{per}$ for $k \in \mathcal{T}_N$, with probability $\eta (1-\eta)^{N^d-1}$. 

In this case $w_{i}^{\eta,N,\tilde{\omega}} = w_{i}^{1,k,N}$ solves the following problem, which we call here a ``one defect'' supercell problem:
\begin{eqnarray} \label{onedef}
\left  \{
\begin{aligned}
& -\mathrm{div}\left(\left(A_{per} + \mathds{1}_{\{Q+ k\}}C_{per}\right)(\nabla w_{i}^{1,k,N} + e_i)\right) = 0\quad \mathrm{in} \, \, I_N,\\
& w_{i}^{1,k,N} \, \, (N\ZZ)^d-\mathrm{periodic}.
\end{aligned}
\right.
\end{eqnarray} 

\item $A_{\eta}(x,\tilde{\omega}) = A_{per} + \mathds{1}_{\{Q+l\} \cup \{Q+m\}}C_{per} \; \mathrm{for} \; (l,m) \in \mathcal{T}_N$, $l \neq m$, with probability \newline $\eta^2 (1-\eta)^{N^d-2}$. 

In this case $w_{i}^{\eta,N,\tilde{\omega}} = w_{i}^{2,l,m,N}$ solves the following problem, which we call here a ``two defects'' supercell problem:
\begin{eqnarray} \label{twodef}
\left  \{
\begin{aligned}
& -\mathrm{div}\left(\left(A_{per} + \mathds{1}_{\{Q+ l\}\cup \{Q+m\}}C_{per}\right)(\nabla w_{i}^{2,l,m,N} + e_i)\right) = 0\quad \mathrm{in} \, \, I_N,\\
& w_{i}^{2,l,m,N} \, \, (N\ZZ)^d-\mathrm{periodic}. 
\end{aligned}
\right.
\end{eqnarray} 
\end{itemize}

All the other possible values for $A_{\eta}$, which are of probability less than $\eta^3$ and which we will not use in this article, can be obtained using similar computations.\\

\begin{figure}[H]
\center
\begin{tabular}{ccc}
\includegraphics[width=4.5cm, height=4.5cm]{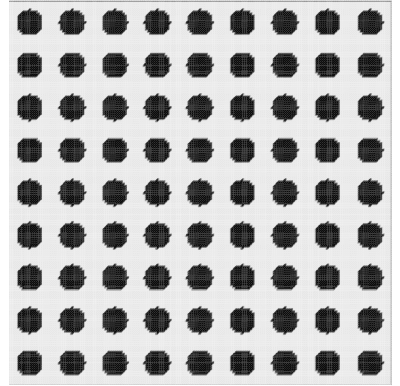} & \includegraphics[width=4.5cm, height=4.5cm]{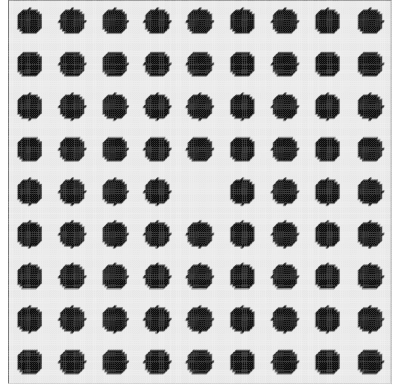} & \includegraphics[width=4.5cm, height=4.5cm]{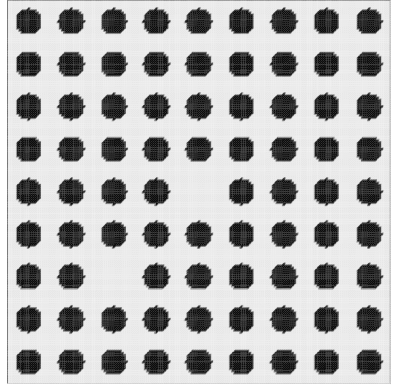} 
\end{tabular}
\caption{From left to right: zero defect, one defect and two defects.}
\label{illus}
\end{figure}

An instance of a setting with zero, one and two defects is shown in Figure \ref{illus} in the two-dimensional case of a material $A_{per}$ consisting of a lattice of inclusions.

Let us define $A_1^k = A_{per} + \mathds{1}_{\{Q+ k\}}C_{per}$ and $A_2^{l,m} = A_{per} + \mathds{1}_{\{Q+ l\}\cup \{Q+m\}}C_{per}$. Then (\ref{expect}) reads
\begin{eqnarray*}
A_{\eta}^*e_i = \lim_{N \rightarrow \infty} \left( \frac{(1-\eta)^{N^d}}{N^d} \int_{I_N}A_{per}(\nabla w_{i}^{0} + e_i ) + \sum_{k \in \mathcal{T}_N}\frac{\eta (1-\eta)^{N^d-1}}{N^d} \int_{I_N}A_1^k(\nabla w_{i}^{1,k,N} + e_i ) \right.  \\
 + \left. \sum_{l,m \in \mathcal{T}_N, l\neq m}\frac{\eta^2 (1-\eta)^{N^d-2}}{2 N^d} \int_{I_N}A_2^{l,m}(\nabla w_{i}^{2,l,m,N} + e_i ) + \cdots \right).
\end{eqnarray*}

It is clear, by $(N\ZZ)^d$-periodicity, that $\displaystyle \int_{I_N}A_1^k(\nabla w_{i}^{1,k,N} + e_i )$ does not depend on the position $k \in \mathcal{T}_N$ of the defect. Likewise,  $\displaystyle \int_{I_N}A_2^{l,m}(\nabla w_{i}^{2,l,m,N} + e_i )$ only depends on the vector $m-l$. Thus we can rewrite
\begin{eqnarray} \label{predevn}
A_{\eta}^*e_i = \lim_{N \rightarrow \infty} \left( (1-\eta)^{N^d}A_{per}^* e_i + \eta (1-\eta)^{N^d-1}\int_{I_N}A_1^0(\nabla w_{i}^{1,0,N} + e_i ) \right. \nonumber \\ + \left.\sum_{k \in \mathcal{T}_N, k\neq 0} \frac{\eta^2 (1-\eta)^{N^d-2}}{2}\int_{I_N}A_2^{0,k}(\nabla w_{i}^{2,0,k,N} + e_i ) + \cdots \right).
\end{eqnarray}

This is of the form
\begin{eqnarray} \label{devn}
\begin{aligned}
A_{\eta}^* &= \lim_{N \rightarrow \infty} \sum_{p=0}^{N^d} \eta^p A_{p}^{*,N}\\
           &= \lim_{N \rightarrow \infty} \left(A_{0}^{*,N} + \eta A_{1}^{*,N} + \eta^2  A_{2}^{*,N} + o_N(\eta^2)\right),
\end{aligned}
\end{eqnarray}
where the remainder $o_N(\eta^2)$ depends on $N$.\\

Explicitly expanding the polynomials in $\eta$ up to the second-order in (\ref{predevn}), we obtain:
\begin{eqnarray} \label{terma0}
 A_{0}^{*,N} = A_{per}^*,
\end{eqnarray}
\begin{eqnarray} \label{terma1}
A_{1}^{*,N} e_i = \int_{I_N}A_1^0(\nabla w_{i}^{1,0,N} + e_i ) - \int_{I_N}A_{per}(\nabla w_{i}^{0} + e_i ),
\end{eqnarray}
\begin{eqnarray} \label{terma2}
\begin{aligned}
A_{2}^{*,N} e_i = \frac{1}{2}\sum_{k \in \mathcal{T}_N, k\neq 0}  \left(\int_{I_N}A_2^{0,k}(\nabla w_{i}^{2,0,k,N} + e_i ) - 2 \int_{I_N}A_1^0(\nabla w_{i}^{1,0,N} + e_i ) \right. \\
 + \left. \int_{I_N}A_{per}(\nabla w_{i}^{0} + e_i ) \right)
\end{aligned}
\end{eqnarray}
as the first three coefficients in (\ref{devn}).

\begin{remark}
The structure of $A_{p}^{*,N}$ for $p \in \mathbb{N}$ is obviously related to that of the polynomial $(1-x)^p$.
\end{remark}

Our approach consists in formally exchanging the limits $N \rightarrow \infty$ and $\eta \rightarrow 0$ in (\ref{devn}). In the next section, we show that $A_{1}^{*,N}$ is a converging sequence when $N \rightarrow \infty$. The case of $A_{2}^{*,N}$, which is shown to be a bounded sequence and thus to converge up to extraction, is discussed in Section \ref{seco}.\\

 We are not able to prove, though, that $A_{\eta}^* - \displaystyle \lim_{N \rightarrow \infty} (A_{per}^* - \eta A_{1}^{*,N} - \eta^2 A_{2}^{*,N}) = o(\eta^2)$ with a remainder term $o(\eta^2)$ independent of $N$.

\begin{remark}
The expression of $A_1^{*,N}$ (and likewise $A_2^{*,N}$) is reminiscent of standard expressions in solid state theory: each of the two integrals in the definition (\ref{terma1}) of $A_{1}^{*,N}$ scales as the volume $N^d$ of the domain $I_N$, and a priori needs to be renormalized in order to give a finite limit. The difference however has a finite limit without renormalization. In solid state physics, it is common to substract a jellium, that is, a uniform background, and proceed similarly.
\end{remark}

\subsection{Convergence of the first-order term $A_1^{*,N}$}

We study here the convergence, as $N$ goes to infinity, of $A_{1}^{*,N}$ defined by (\ref{terma1}),
and prove:
\begin{Prop} \label{conva1}
$A_{1}^{*,N}$ converges to a finite limit $\bar{A_{1}^{*}}$ in $\RR^{d \times d}$ when $N \rightarrow \infty$.
\end{Prop}

\begin{proof}

We fix $(i,j)$ in $\llbracket 1, d \rrbracket^2$ and study the convergence of $A_{1}^{*,N} e_i \cdot e_j$.\\

Let us define the adjoint problems to the cell problems (\ref{cellper}):
\begin{eqnarray} \label{celladj}
\left  \{
\begin{aligned}
& -\mathrm{div}\left(A_{per}^T(\nabla \tilde{w}_{j}^0 + e_j ) \right) = 0 \quad \mathrm{in} \, \, Q,\\
& \tilde{w}_{j}^0 \; \; \mathbb{Z}^d-\mathrm{periodic}, 
\end{aligned}
\right.
\end{eqnarray}
where we have denoted by $A_{per}^T$ the transposed matrix of $A_{per}$. Then using (\ref{onedef}) and the definition of $A_1^0$, we have
\begin{eqnarray*}
\int_{I_N}A_1^0(\nabla w_{i}^{1,0,N} + e_i ) \cdot e_j &=&  \int_{I_N}A_1^0(\nabla w_{i}^{1,0,N} + e_i ) \cdot (e_j + \nabla \tilde{w}_{j}^0) \\
&=& \int_{I_N}A_{per}(\nabla w_{i}^{1,0,N} + e_i ) \cdot (e_j + \nabla \tilde{w}_{j}^0)  \\
& & + \int_{Q}C_{per}(\nabla w_{i}^{1,0,N} + e_i ) \cdot (e_j + \nabla \tilde{w}_{j}^0).
\end{eqnarray*}

Next, using (\ref{celladj}), we note that
\begin{eqnarray*}
\int_{I_N}A_{per}(\nabla w_{i}^{1,0,N} + e_i ) \cdot (e_j + \nabla \tilde{w}_{j}^0) &=&  \int_{I_N}(\nabla w_{i}^{1,0,N} + e_i ) \cdot A_{per}^T (e_j + \nabla \tilde{w}_{j}^0) \\
&=& \int_{I_N} e_i \cdot A_{per}^T (e_j + \nabla \tilde{w}_{j}^0) \\
&=& N^d \left(A_{per}^T\right)^* e_j \cdot e_i, 
\end{eqnarray*}
and applying (\ref{homtens}) to the periodic tensor $A_{per}^T$ and noticing that $(A_{per}^T)^* = (A_{per}^*)^T$, we obtain

\begin{eqnarray} \label{a1nd}
\int_{I_N}A_1^0(\nabla w_{i}^{1,0,N} + e_i ) \cdot e_j = N^d A_{per}^* e_i \cdot e_j + \int_{Q}C_{per}(\nabla w_{i}^{1,0,N} + e_i ) \cdot (e_j + \nabla \tilde{w}_{j}^0).
\end{eqnarray}

Since, by definition, 
\begin{eqnarray*} 
A_{1}^{*,N} e_i &=& \int_{I_N}A_1^0(\nabla w_{i}^{1,0,N} + e_i ) - \int_{I_N}A_{per}(\nabla w_{i}^{0} + e_i)\\
&=& \int_{I_N}A_1^0(\nabla w_{i}^{1,0,N} + e_i ) - N^d A_{per}^* e_i ,
\end{eqnarray*} 
we deduce from (\ref{a1nd}) that
\begin{eqnarray} \label{adef1}
A_{1}^{*,N} e_i \cdot e_j = \int_{Q}C_{per}(\nabla w_{i}^{1,0,N} + e_i ) \cdot (e_j + \nabla \tilde{w}_{j}^0).
\end{eqnarray}

We now define 
\begin{eqnarray} \label{compare1}
q_i^{1,0,N} = w_{i}^{1,0,N} - w_i^0,
\end{eqnarray}
which solves
\begin{eqnarray} \label{perturb}
\left  \{
\begin{aligned}
& -\mathrm{div} \left(A_1^0 \nabla q_i^{1,0,N} \right) = \mathrm{div} (\mathds{1}_Q C_{per} (\nabla w_i^0 + e_i)) \quad \mathrm{in} \; I_N,\\
& q_i^{1,0,N}\,(N\ZZ)^d-\mathrm{periodic}. 
\end{aligned}
\right.
\end{eqnarray} 
We deduce from Lemma \ref{conv1} of the appendix, applied to (\ref{perturb}), that $ \nabla q_i^{1,0,N}$ converges in $L^2_{loc}(\mathbb{R}^d)$, when $N \rightarrow + \infty$, to $\nabla q_i^{1,0,\infty}$, where $q_i^{1,0,\infty}$ is a $L^2_{loc}(\RR^d)$ function solving
\begin{eqnarray} \label{perturbdef} 
\left  \{
\begin{aligned}
& -\mathrm{div} \left(A_1^0 \nabla q_i^{1,0,\infty} \right) = \mathrm{div} (\mathds{1}_Q C_{per} (\nabla w_i^0 + e_i)) \quad \mathrm{in} \; \RR^d,\\
& \nabla q_i^{1,0,\infty} \in L^2(\RR^d). 
\end{aligned}
\right.
\end{eqnarray} 

Defining $w_{i}^{1,0,\infty} = w_i^0 + q_i^{1,0,\infty}$, it is clear that $\nabla w_{i}^{1,0,N}$ converges in $L^2(Q)$ to $\nabla w_{i}^{1,0,\infty}$. It follows from (\ref{adef1}) that $A_{1}^{*,N} \underset{N \rightarrow + \infty}{\rightarrow} \bar{A_{1}^{*}}$ with $\bar{A_{1}^{*}}$ defined by 
\begin{eqnarray} \label{Ainfiny}
\forall(i,j) \in \llbracket 1,d \rrbracket^2, \; \; \bar{A_{1}^{*}} e_i \cdot e_j &=& \int_{Q}C_{per}(\nabla w_{i}^{1,0,\infty} + e_i ) \cdot (e_j + \nabla \tilde{w}_{j}^0).
\end{eqnarray}
\end{proof}

The computation of $\bar{A_1^{*}}$ requires to solve (\ref{perturbdef}) which is defined in $\RR^d$, but, in sharp contrast to the stochastic cell problems (\ref{etacell}), is deterministic and has a right-hand side with compact support in $\RR^d$. In practice, problem (\ref{perturbdef}) is truncated on $I_N$. The following result gives insight on the truncation error.  

\begin{lemma} \label{rate}
Assume that $d \geq 3$ and that the unit cell $Q$ contains an inclusion $D$, the boundary of which has regularity $\mathcal{C}^{1,\mu}$ for some $0 < \mu < 1$, and such that $\mathrm{dist}(D,\partial Q)>0$. Assume also that $A_{per}$ is Hölder continuous in $\overline{D}$ and in $\overline{Q \backslash D}$. Then there exists a tensor $B_1^{*,N}$, computed on $I_N$, and a constant $K$ independent of $N$ such that 
$$ |B_1^{*,N} - \bar{A_{1}^{*}}| \leq K N^{-d}.$$
\end{lemma}

\begin{proof}

\textbf{Step 1.}\\

Fix $(i,j)$ in $\llbracket 1, d \rrbracket^2$. We first define the adjoint problem for (\ref{perturb}), namely
\begin{eqnarray} \label{qadjoint}
\left  \{
\begin{aligned}
& -\mathrm{div} \left((A_1^0)^T \nabla \tilde{q}_j^{1,0,N}\right) = \mathrm{div} (\mathds{1}_Q C_{per}^T (\nabla \tilde{w}_j^0 + e_j)) \quad \mathrm{in} \; I_N,\\
& \tilde{q}_j^{1,0,N}\,(N \ZZ)^d-\mathrm{periodic}.
\end{aligned}
\right.
\end{eqnarray} 

Applying Lemma \ref{conv1} to (\ref{qadjoint}), we also introduce the limit $\tilde{q}_{j}^{1,0,\infty}$ of $\tilde{q}_{j}^{1,0,N}$ when $N \rightarrow \infty$. It solves the adjoint problem of (\ref{perturbdef}).\\

 Then, using (\ref{qadjoint}), we obtain
\begin{eqnarray*}
\int_{Q}C_{per}\nabla q_{i}^{1,0,N} \cdot (e_j + \nabla \tilde{w}_{j}^0)&=&\int_{I_N}\nabla q_{i}^{1,0,N} \cdot \mathds{1}_Q C_{per}^T  (e_j + \nabla \tilde{w}_{j}^0)\\
&=& - \int_{I_N} \nabla q_{i}^{1,0,N} \cdot (A_1^0)^T \nabla \tilde{q}_j^{1,0,N} \\
&=& - \int_{I_N} A_1^0 \nabla q_{i}^{1,0,N} \cdot \nabla \tilde{q}_j^{1,0,N}.
\end{eqnarray*}   
Consequently, (\ref{adef1}) and the definition (\ref{compare1}) of $q_i^{1,0,N}$ yield
\begin{eqnarray} \label{rwa1}
A_{1}^{*,N} e_i \cdot e_j &=& \int_{Q}C_{per}(\nabla w_{i}^{0} + e_i ) \cdot (e_j + \nabla \tilde{w}_{j}^0) - \int_{I_N} A_1^0 \nabla q_{i}^{1,0,N} \cdot \nabla \tilde{q}_j^{1,0,N}.
\end{eqnarray}
We  know from Lemma \ref{conv1} applied to (\ref{perturb}) and (\ref{qadjoint}) that $\mathds{1}_{I_N} \nabla q_{i}^{1,0,N}$  and $\mathds{1}_{I_N} \nabla \tilde{q}_j^{1,0,N}$converge strongly in $L^2(\RR^d)$ to $\nabla q_{i}^{1,0,\infty}$ and $\nabla \tilde{q}_j^{1,0,\infty}$ respectively, when $N \rightarrow \infty$. Passing to the limit in (\ref{rwa1}) then gives 
\begin{eqnarray*}
\bar{A_{1}^{*}} e_i \cdot e_j &=& \int_{Q}C_{per}(\nabla w_{i}^{0} + e_i ) \cdot (e_j + \nabla \tilde{w}_{j}^0) - \int_{\RR^d} A_1^0 \nabla q_{i}^{1,0,\infty} \cdot \nabla \tilde{q}_j^{1,0,\infty}.
\end{eqnarray*}

We now define $v_{i}^{1,0,N}$ and $\tilde{v}_j^{1,0,N}$ solutions to (\ref{perturb}) and (\ref{qadjoint}) with homogeneous Dirichlet (instead of periodic) boundary conditions on the boundary $\partial I_N$ of $I_N$, and the tensor~$B_1^{*,N}$ by
\begin{eqnarray} \label{B}
B_{1}^{*,N} e_i \cdot e_j &=& \int_{Q}C_{per}(\nabla w_{i}^{0} + e_i ) \cdot (e_j + \nabla \tilde{w}_{j}^0) - \int_{I_N} A_1^0 \nabla v_{i}^{1,0,N} \cdot \nabla \tilde{v}_j^{1,0,N}.
\end{eqnarray}

The proof of Proposition \ref{conva1} is easily adapted to show that $B_{1}^{*,N}$ converges to $A_{1}^{*}$ as $N$ goes to infinity. \\

\textbf{Step 2.\\}

 We consider
\begin{eqnarray*}
\left(B_{1}^{*,N} - \bar{A_{1}^{*}}\right)e_i \cdot e_j &=& \int_{\RR^d} A_1^0 \nabla q_{i}^{1,0,\infty} \cdot \nabla \tilde{q}_j^{1,0,\infty} - \int_{I_N} A_1^0 \nabla v_{i}^{1,0,N} \cdot \nabla \tilde{v}_j^{1,0,N},   
\end{eqnarray*}
and expand the difference $B_{1}^{*,N} - \bar{A_{1}^{*}}$ as follows:
\begin{eqnarray} \label{cut}
\begin{aligned}
\hspace{-0.5cm}\left(B_{1}^{*,N} - \bar{A_{1}^{*}}\right)e_i \cdot e_j = \int_{\RR^d \backslash I_N} A_1^0 \nabla q_{i}^{1,0,\infty} \cdot \nabla \tilde{q}_j^{1,0,\infty}&  \\ 
 &\hspace{-5cm}+ \left(\int_{I_N} A_1^0 \nabla q_{i}^{1,0,\infty} \cdot \nabla \tilde{q}_j^{1,0,\infty} - \int_{I_N} A_1^0 \nabla v_{i}^{1,0,N} \cdot \nabla \tilde{v}_j^{1,0,N} \right).
\end{aligned}
\end{eqnarray}

We now show that the two terms in the right-hand side of (\ref{cut}) converge to $0$ as $N^{-d}$ when $N \rightarrow +\infty$.\\

We first note that the results of Lemma \ref{regdec} of the appendix, stated for a $\ZZ^d$-periodic matrix, can be readily extended to address $A_1^0$ since $A_1^0$ is equal to $A_{per}$ in $\RR^d \backslash Q$.\\

We deduce from Lemma \ref{laxlem} applied to (\ref{perturbdef}) that $q_i^{1,0,\infty}$ is defined uniquely up to an additive constant. Moreover, $A_{per}$ being piecewise Hölder continuous, we deduce from Lemma \ref{regdec} that there exists a unique solution to (\ref{perturbdef}) which converges to zero at infinity.\\

 Since we only use $\nabla q_i^{1,0,\infty}$ in $\bar{A_{1}^{*}}$, we can thus assume without loss of generality that $q_i^{1,0,\infty}$ converges to zero at infinity. Likewise, we assume that 
$\tilde{q}_j^{1,0,\infty}$ converges to zero at infinity.\\

We then deduce from Lemma \ref{regdec} that there exists a constant $K$ independent of $N$ such that for $|x| \geq 1$,
\begin{eqnarray}
&|q_i^{1,0,\infty}(x)| \leq K|x|^{1-d}, \; |\tilde{q}_j^{1,0,\infty}(x)| \leq K|x|^{1-d}, \label{qbound} \\
&|\nabla q_i^{1,0,\infty}(x)| \leq K|x|^{-d}, \; |\nabla \tilde{q}_j^{1,0,\infty}(x)| \leq K|x|^{-d}, \label{nablaqbound} \\
 & |v_i^{1,0,N}(x)| \leq K |x|^{1-d}, \; |\tilde{v}_j^{1,0,N}(x)| \leq K |x|^{1-d} \label{vbound}, \\
& |\nabla v_i^{1,0,N}(x)| \leq K|x|^{-d}, \;|\nabla \tilde{v}_j^{1,0,N}(x)| \leq K|x|^{-d} \label{nablavbound}. 
\end{eqnarray}

Using (\ref{nablaqbound}), we have $\|\nabla q_{i}^{1,0,\infty}\|_{L^2(\RR^d \backslash I_N)}\leq K N^{-d/2}$ and $\|\nabla \tilde{q}_j^{1,0,\infty} \|_{L^2(\RR^d \backslash I_N)} \leq K N^{-d/2}$, and
\begin{eqnarray}\label{term1}
\begin{aligned}
\left| \int_{\RR^d \backslash I_N} A_1^0 \nabla q_{i}^{1,0,\infty} \cdot \nabla \tilde{q}_j^{1,0,\infty}  \right| &\leq \beta \|\nabla q_{i}^{1,0,\infty}\|_{L^2(\RR^d \backslash I_N)}\|\nabla \tilde{q}_j^{1,0,\infty} \|_{L^2(\RR^d \backslash I_N)} \\
&\leq KN^{-d},
\end{aligned}
\end{eqnarray}
where $\beta$ is defined in (\ref{betabound}).\\

We now address the second term of the right-hand side of (\ref{cut}) and write
\begin{eqnarray*}
\int_{I_N} A_1^0 \nabla v_{i}^{1,0,N} \cdot \nabla \tilde{v}_j^{1,0,N} -  \int_{I_N} A_1^0 \nabla q_{i}^{1,0,\infty} \cdot \nabla \tilde{q}_j^{1,0,\infty} & & \\
&& \hspace{-8cm} = \int_{I_N}A_1^0 (\nabla v_{i}^{1,0,N} - \nabla q_{i}^{1,0,\infty}) \cdot  \nabla \tilde{v}_j^{1,0,N} + \int_{I_N} A_1^0 \nabla q_{i}^{1,0,\infty} \cdot (\nabla \tilde{v}_j^{1,0,N}-\nabla \tilde{q}_{j}^{1,0,\infty}).
\end{eqnarray*}
Since $\mathrm{div}\left(A_1^0 (\nabla v_{i}^{1,0,N} - \nabla q_{i}^{1,0,\infty})\right) = \mathrm{div} \left((A_1^0)^T (\nabla \tilde{v}_{j}^{1,0,N} - \nabla \tilde{q}_{j}^{1,0,\infty})\right)=0$ in $I_N$, and $\tilde{v}_{j}^{1,0,N} = 0$ on $\partial I_N$, we have, using integration by parts,
\begin{eqnarray*}
\int_{I_N}A_1^0 (\nabla v_{i}^{1,0,N} - \nabla q_{i}^{1,0,\infty}) \cdot  \nabla \tilde{v}_j^{1,0,N} + \int_{I_N} A_1^0 \nabla q_{i}^{1,0,\infty} \cdot (\nabla \tilde{v}_j^{1,0,N}-\nabla \tilde{q}_{j}^{1,0,\infty}) & & \\
&& \hspace{-9cm} =  \int_{\partial I_N} A_1^0 (\nabla v_{i}^{1,0,N} - \nabla q_{i}^{1,0,\infty}) \cdot \nu \, \tilde{v}_j^{1,0,N}\\
&& \hspace{-8cm} + \int_{\partial I_N} (A_1^0)^T (\nabla \tilde{v}_{j}^{1,0,N} - \nabla \tilde{q}_{j}^{1,0,\infty}) \cdot \nu \, q_{i}^{1,0,\infty} \\
&& \hspace{-9cm} =  \int_{\partial I_N} (A_1^0)^T (\nabla \tilde{v}_{j}^{1,0,N} - \nabla \tilde{q}_{j}^{1,0,\infty}) \cdot \nu \, q_{i}^{1,0,\infty}, 
\end{eqnarray*}
where $\nu$ is the unit outward normal vector to $\partial I_N$.\\

The estimates (\ref{qbound}) and (\ref{nablavbound}) imply
$$ \|q_{i}^{1,0,\infty}\|_{L^{\infty}(\partial I_N)} \leq K N^{1-d}, \quad \|(A_1^0)^T (\nabla \tilde{v}_{j}^{1,0,N} - \nabla \tilde{q}_{j}^{1,0,\infty}) \cdot \nu\|_{L^{\infty}(\partial I_N)} \leq K N^{-d},$$
while the measure of the boundary $\partial I_N$ scales as $N^{1-d}$. Hence
$$ \left|\int_{\partial I_N} (A_1^0)^T (\nabla \tilde{v}_{j}^{1,0,N} - \nabla \tilde{q}_{j}^{1,0,\infty}) \cdot \nu \, q_{i}^{1,0,\infty} \right| \leq K N^{-d},$$
and then 
\begin{eqnarray}\label{term2}
\left|\int_{I_N} A_1^0 \nabla v_{i}^{1,0,N} \cdot \nabla \tilde{v}_j^{1,0,N} -  \int_{I_N} A_1^0 \nabla q_{i}^{1,0,\infty} \cdot \nabla \tilde{q}_j^{1,0,\infty} \right| \leq K N^{-d}.
\end{eqnarray}

We conclude by substituting (\ref{term1}) and (\ref{term2}) into (\ref{cut}).

\end{proof}

\begin{remark}
We assume $d \geq 3$ and piecewise Hölder regularity on $A_{per}$, and use Dirichlet boundary conditions in Lemma \ref{rate}, because our proof relies on Lemma \ref{regdec}. Note however that the numerical experiments of Section \ref{num} show, in dimension $d=2$, that we again obtain the rate $N^{-d}$ in the convergence of $A_1^{*,N}$ to $\bar{A_1^{*}}$ for two different $A_{per}$, one being piecewise Hölder continuous in the sense of Lemma \ref{rate} and the other not, and with periodic boundary conditions. Moreover, the explicit computations of Proposition \ref{preuve1d} show that in dimension one, and without any assumption of regularity on $A_1^0$, the rate of convergence of $A_1^{*,N}$ to $\bar{A_1^{*}}$ is $N^{-1}$.
\end{remark}

\subsection{Second-order term} \label{seco}

For completeness, we state here the result regarding $A_{2}^{*,N}$ proved in \cite{these_Arnaud}:
\begin{Prop} \label{conva2}
The sequence $A_{2}^{*,N}$ defined by (\ref{terma2}) is bounded in $\RR^{d\times d}$ and therefore converges up to extraction.
\end{Prop}

We strongly believe that $A_{2}^{*,N}$ is actually a convergent sequence, as shown by our numerical tests thereafter. In fact, we even believe (see \cite{these_Arnaud}) that we can write the expression of the limit. Note also that the explicit computations of Section \ref{1dc} prove the convergence of $A_{2}^{*,N}$ in dimension one.

\section{Numerical experiments}\label{num}

Our purpose in this section is to assess the approximation of $A_{\eta}^*$ by the second-order expansion $A_{per}^{*} + \eta A_1^{*,N} + \eta^2 A_2^{*,N}$. The limited computational facilities we have access to impose that we restrict ourselves to the two-dimensional case. We first explain our general methodology and then make precise the specific settings.

\subsection{Methodology}

We will consider two commonly used composite materials as periodic reference materials~$A_{per}$. The first material consists of a constant background reinforced by a periodic lattice of circular inclusions, that is
$$ A_{per}(x_1,x_2) = 20 \times Id + 100 \sum_{k \in \ZZ^2} \mathds{1}_{B(k,0.3)}(x_1,x_2) \times Id,$$
where $B(k,0.3)$ is the ball of center $k$ and radius $0.3$. The second material is a laminate for which
$$ A_{per}(x_1,x_2) = 20 \times Id + 100 \sum_{l \in \ZZ} \mathds{1}_{l \leq x_1 \leq l+1}(x_1,x_2) \times Id.$$ 
In the case of material $1$, the role of the perturbation is, loosely speaking, to randomly eliminate some fibers:
$$ C_{per}(x_1,x_2) = -100 \sum_{k \in \ZZ^d} \mathds{1}_{B(k,0.3)}(x_1,x_2) \times Id.$$
In the case of material $2$, the perturbation consists in a random modification of the lamination direction:
$$ C_{per}(x_1,x_2) = - 100 \sum_{l \in \ZZ} \mathds{1}_{l \leq x_1 \leq l+1}(x_1,x_2) \times Id + 100 \sum_{l \in \ZZ} \mathds{1}_{l \leq x_2 \leq l+1}(x_1,x_2) \times Id.$$

In both cases, we have chosen the coefficients $20$ and $100$ in order to have a high contrast between $A_{per}$ and $A_{per} + C_{per}$, and thus for the perturbation to be significant. There is of course nothing specific in the actual value of these coefficients.\\

\begin{figure}[h] 
\center
\begin{tabular}{cc}
\includegraphics[width=0.42\textwidth]{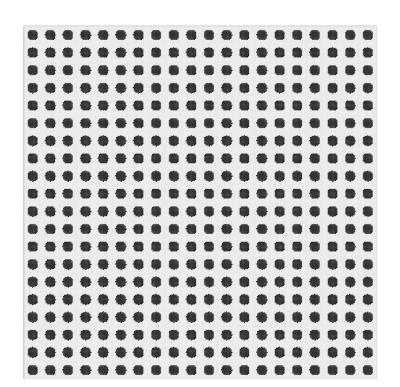} & \includegraphics[width=0.42\textwidth]{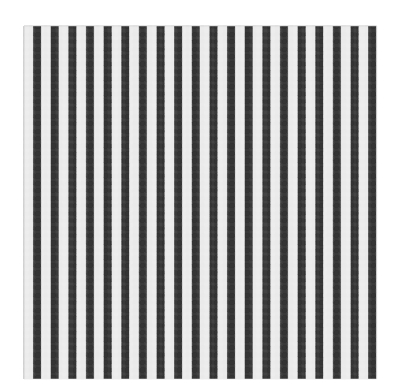}
\end{tabular}
\caption{Left: a periodic lattice of circular inclusions. Right: a one-dimensional laminate.}
\label{fibres}
\end{figure}

These two materials are shown in Figure \ref{fibres}.\\ 

Our goal is to compare $A_{\eta}^*$ with its approximation $A_{per}^{*} + \eta A_1^{*,N} + \eta^2 A_2^{*, N}$ for each of these two particular settings. A major computational difficulty is the computation of the ``exact'' matrix $A_{\eta}^*$ given by formula (\ref{homogal}). It ideally requires to solve the stochastic cell problems (\ref{etacell}) on $\RR^d$. To this end we first use ergodicity and formulae (\ref{ergodiccorr}) and (\ref{expect}), and actually compute, for a given realization $\omega$ and a domain $I_N$ which is here equal to $[0,N]^2$ for convenience, $A_{\eta}^{*,N}(\omega)$ defined by 
\begin{eqnarray} \label{ergodicint} 
A_{\eta}^{*,N}(\omega) e_i = \frac{1}{N^d} \int_{I_N}A_{\eta}(x,\omega)(\nabla w_{i}^{\eta,N,\omega}(x) + e_i ) dx. 
\end{eqnarray} 

In a second step, we take averages over the realizations $\omega$.\\

For each $\omega$, we use the finite element software FreeFem++ (available at www.freefem.org) to solve the boundary value problems (\ref{ergodiccorr}) and compute the integrals (\ref{ergodicint}). We work with standard P1 finite elements on a triangular mesh such that there are $10$ degrees of freedom on each edge of the unit cell $Q$.\\  

We define an approximate value $A_{\eta}^{*,N}$ as the average of $A_{\eta}^{*,N}(\omega)$ over $40$ realizations~$\omega$. Our numerical experiments indeed show that the number $40$ is sufficiently large for the convergence of the Monte-Carlo computation. We then let $N$ grow from $5$ to $80$ by increments of $5$. We observe that $A_{\eta}^{*,N}$ stabilizes at a fixed value around $N=80$ and thus take $A_{\eta}^{*,80}$ as the reference value for $A_{\eta}^{*}$ in our subsequent tests.\\

The next step is to compute the zero-order term $A_{per}^*$, and the first-order and second-order deterministic corrections $A_{1}^{*,N}$ and $A_{2}^{*,N}$. Using the same mesh and finite elements as for our reference computation above, we compute $A_{per}^*$ using (\ref{cellper}) and (\ref{homtens}), and for each~$N$ we compute $A_1^{*,N}$ and $A_2^{*,N}$ using (\ref{terma1}) and (\ref{terma2}). We again let $N$ grow from $5$ to $80$ by increments of $5$ for $A_1^{*,N}$. The computation of $A_2^{*,N}$ being significantly more expensive (note that in (\ref{terma2}) there is not only an integral over $I_N$ but also a sum over the $N^2$ cells) we have to limit ourselves to $N=25$ and approximate the value for $N$ larger than $25$ by the value obtained for $N=25$.\\

Before presenting our results, we wish to discuss our expectations. Note that there are three distinct sources of error:
\begin{itemize}
\item the finite elements discretization error; 
\item the truncation error due to the replacement of $\RR^d$ with $I_N$, in the computation of the stochastic cell problems (\ref{etacell}) that are replaced with (\ref{ergodiccorr}), as well as in the computation of the integrals (\ref{ergodicint}); 
\item the stochastic error arising from the approximation of the expectation value (\ref{expect}) by an empirical mean. 
\end{itemize}

The discretization error originates from the fact that, in practice, we only have access to the finite element approximations of all the functions manipulated here (such as $w_i^0$, $w_{i}^{\eta,N,\omega}$,...). Although we have not proved it in the specific context of our work, we believe, because it is shown in a similar weakly random setting (see \cite{Costaouec}), that all the convergences stated here in the infinite-dimensional setting still hold true for the finite-dimensional approximations of the objects. Our numerical results indeed confirm it is the case. In order to eliminate the discretization error from the picture, our practical approach consists in adopting the \emph{same} finite element space for all approximations of the cell and supercell problems, independently of $N$.\\ 

The truncation error is a different issue. For the ``exact'' computation of $A_{\eta}^{*}$ (we mean not using the second-order expansion (\ref{devn}), but (\ref{ergodicint})), we use an empirical mean and a truncation. We know from \cite{BP}, for a \emph{continuous} notion of stationarity analogous to the \emph{discrete} notion (\ref{statio}) we use here, and under mixing conditions which are satisfied in our setting, that the convergence of the truncated approximation to the ideal value holds at a rate $N^{-\kappa}$ with $\kappa$ a non explicit function of the dimension, the mixing exponent and the coercivity constant of the material. On the other hand, in the second-order expansion (\ref{devn}), the zero-order term $A_{per}^*$ is of course free of any truncation error. All that we know for the approximation $A_1^{*,N}$ defined by (\ref{terma1}) to the first-order correction $\bar{A_1^*}$, is stated in Lemma~\ref{rate} in dimension $d \geq 3$, under Hölder regularity assumptions on $A_{per}$, and with Dirichlet boundary conditions replacing periodic ones. One of the aims of our experiments is therefore to draw some numerical conclusions on the convergence of this term when these assumptions are not satisfied. Note that the matrices involved in our test materials are clearly discontinuous functions of $x$. The matrix corresponding to material 1 is piecewise Hölder continuous in the sense of Lemma \ref{rate}, while the matrix corresponding to material 2 is not. As for the second-order approximation $A_2^{*,N}$, we have no insight on the truncation error and we also wish to study its convergence from a numerical point of view.\\

Finally, we have a practical approach to the stochastic error: besides the empirical mean, we provide, for each $N$, the minimum and the maximum values of $A_{\eta}^{*,N}(\omega)$ achieved over the $40$ computations.\\

We now would like to emphasize that the purpose of our numerical tests is not to \emph{prove} that    
$$ A_{\eta}^* = A_{per}^* + \eta A_1^{*,N} + \eta^2 A_2^{*,N} + o(\eta^2)$$
for a remainder term $o(\eta^2)$ that is independent of $N$, of the number of realizations and of the size of the mesh. Establishing experimentally that such an asymptotic holds is too demanding a task. It would indeed require letting $\eta$ go to $0$, which in turn, since we have to observe at least one (and in fact many) event per domain considered, would necessitate a supercell of size $N$ extremely large. We cannot afford such a computational workload.\\

Using our numerical tests, we only hope here to demonstrate, and we indeed do so, that the second-order expansion is an approximation to $A_{\eta}^*$ sufficiently good for all practical purposes, and in particular for $\eta$ not too small ! We will observe that $A_2^{*,N}$ is not only bounded as stated in Proposition \ref{conva2} but, as $N$ goes to infinity, converges to a limit $\bar{A_2^{*}}$, and that both $A_1^{*,N}$ and $A_2^{*,N}$ converge to their respective limits faster than $A_{\eta}^{*,N}$ to $A_{\eta}^{*}$ (which is intuitively expected since the former quantities are deterministic and contain less information). We will also observe that $A_{per}^* + \eta A_1^{*,N}$ is significantly closer to $A_{\eta}^{*}$ than $A_{per}^{*}$, thereby motivating the expansion. The inclusion of the second-order term further improves the situation.

\subsection{Results}

\begin{figure}[h]
\center
\begin{tabular}{cc}
\includegraphics[width=0.42\textwidth]{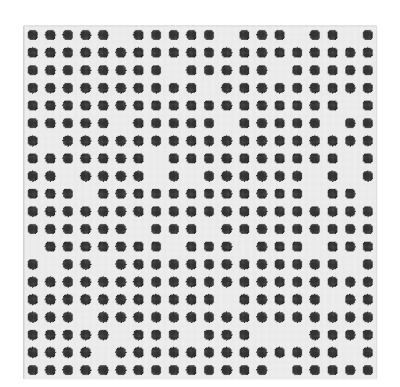} & \includegraphics[width=0.42\textwidth]{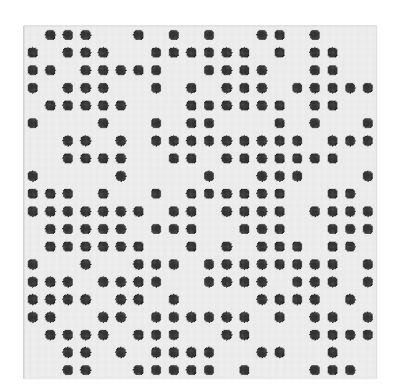} 
\end{tabular}
\caption{Two instances of material $1$ with $\eta=0.1$ (left) and $\eta=0.4$ (right).}
\label{realfibres}
\end{figure}

In order to give an idea on how the perturbation affects the materials considered, we first show some typical realizations in Figure \ref{realfibres} and Figure \ref{reallamine}. Our results are presented in Section \ref{mat1} and Section \ref{mat2} below. Since these results are qualitatively similar for the two materials, we comment on the results altogether in Section \ref{comments}.

\begin{figure}[H]
\center
\begin{tabular}{cc}
\includegraphics[width=0.42\textwidth]{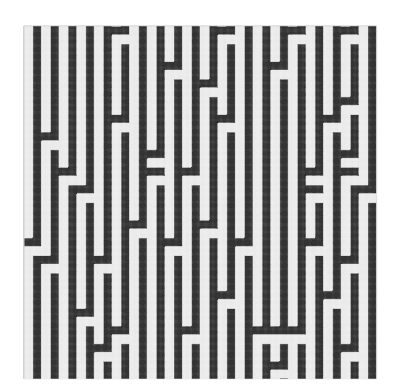} & \includegraphics[width=0.42\textwidth]{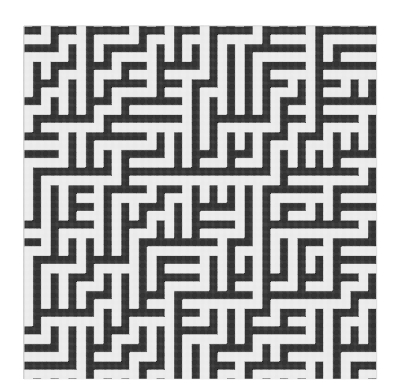} 
\end{tabular}
\caption{Two instances of material $2$ with $\eta=0.1$ (left) and $\eta=0.5$ (right).}
\label{reallamine}
\end{figure}

To present our numerical results, we choose the first diagonal entry $(1,1)$ of all the matrices considered. Other coefficients in the matrices behave qualitatively similarly. As mentioned in the previous Section, we illustrate a practical interval of confidence for our Monte-Carlo computation of $A_{\eta}^*$ by showing, for each $N$, the minimum and maximum values of $A_{\eta}^{*,N}(\omega)$ achieved over the $40$ realizations $\omega$.\\

We will use the following legend in the graphs:
\begin{itemize}
\item {\it periodic}: gives the value of the periodic homogenized tensor $A_{per}^*$;
\item {\it first-order}: gives the value of $A_{per}^* + \eta A_1^{*,N}$;
\item {\it second-order}: gives the value of $A_{per}^* + \eta A_1^{*,N} + \eta^2 A_2^{*,N}$;
\item {\it stochastic mean, minima and maxima}: respectively give the values of $A_{\eta}^{*,N}$ and the extrema obtained in the computation of the empirical mean.
\end{itemize}

Finally, the results are given for \emph{some} specific values of $\eta$ (not necessarily the same for both materials) which serve the purpose of testing our approach in a diversity of situations, from a ``small'' to a ``not so small'' perturbation. 

\subsubsection{Results for material $1$} \label{mat1}

\begin{figure}[H]
\center
\includegraphics[width=14cm, height=10cm]{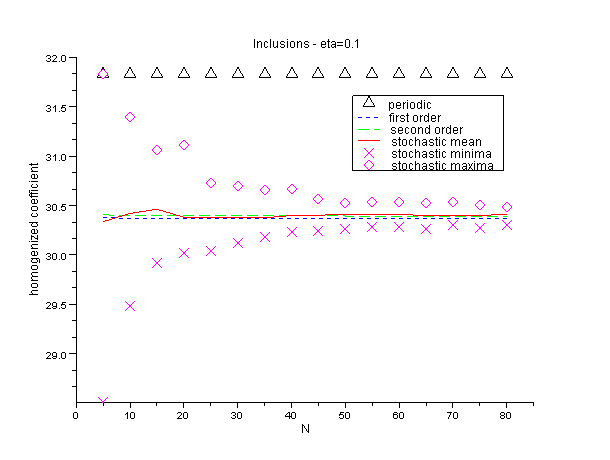} 
\includegraphics[width=14cm, height=10cm]{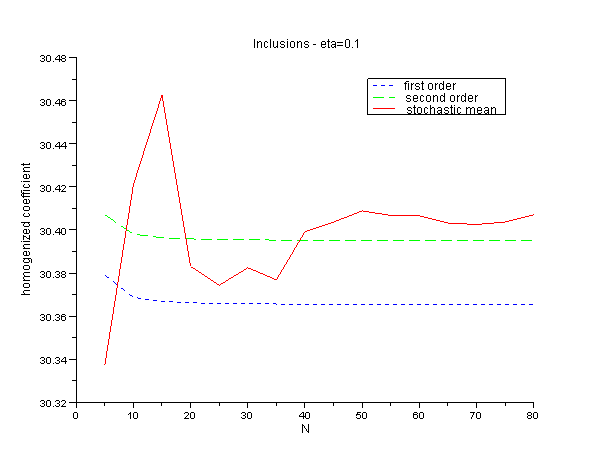} 
\caption{Results for material 1 and $\eta=0.1$. Above: complete results. Below: close-up on $A_{\eta}^{*,N}$ and the first and second-order corrections.}
\label{convfibres01}
\end{figure}

\begin{figure}[H]
\center
\includegraphics[width=14cm, height=11cm]{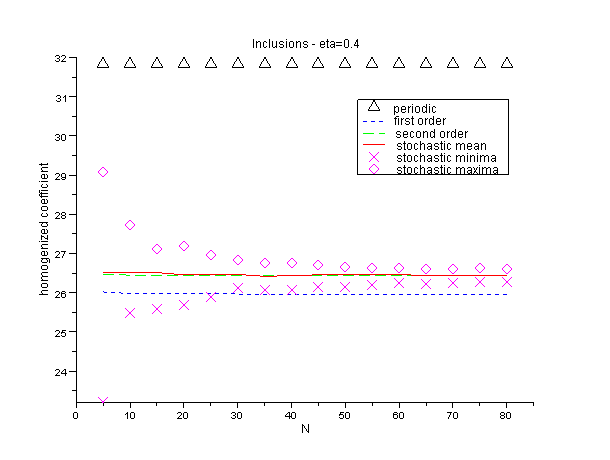} 
\includegraphics[width=14cm, height=11cm]{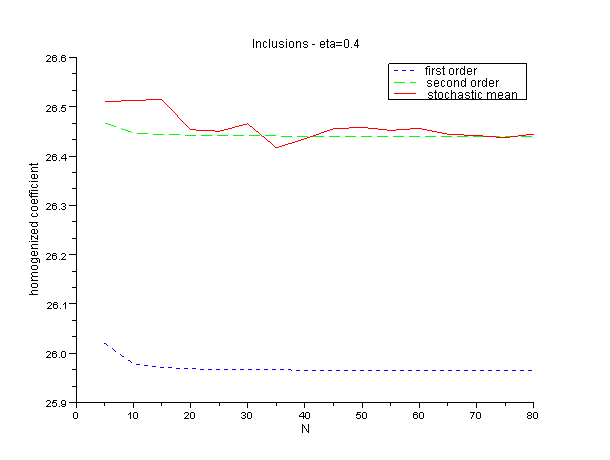} 
\caption{Results for material 1 and $\eta=0.4$. Above: complete results. Below: close-up on $A_{\eta}^{*,N}$ and the first and second-order corrections.}
\label{convfibres04}
\end{figure}

\begin{figure}[H]
\center
\includegraphics[width=14cm, height=11cm]{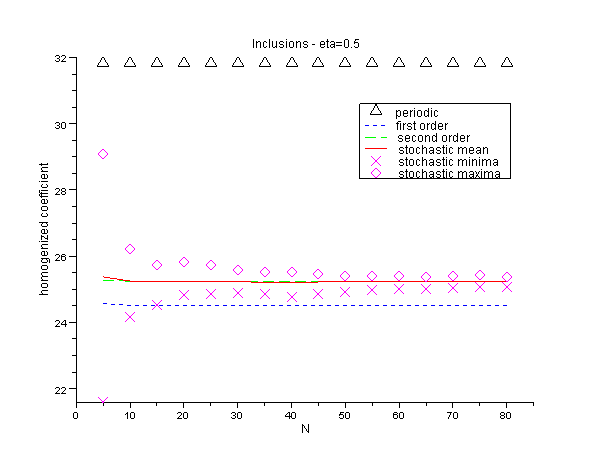} 
\includegraphics[width=14cm, height=11cm]{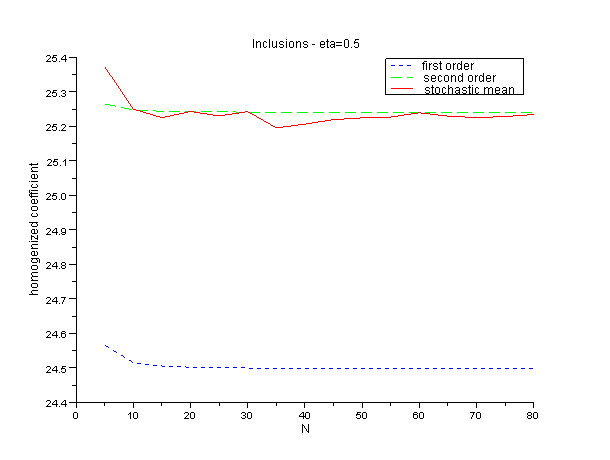} 
\caption{Results for material 1 and $\eta=0.5$. Above: complete results. Below: close-up on $A_{\eta}^{*,N}$ and the first and second-order corrections.}
\label{convfibres05}
\end{figure}

We show the results for $\eta=0.1$, $\eta=0.4$ and $\eta=0.5$ (Figures \ref{convfibres01}, \ref{convfibres04} and \ref{convfibres05} respectively).

\newpage

\newpage

\newpage

\subsubsection{Results for material $2$} \label{mat2}

\begin{figure}[H]
\center
\includegraphics[width=14cm, height=10cm]{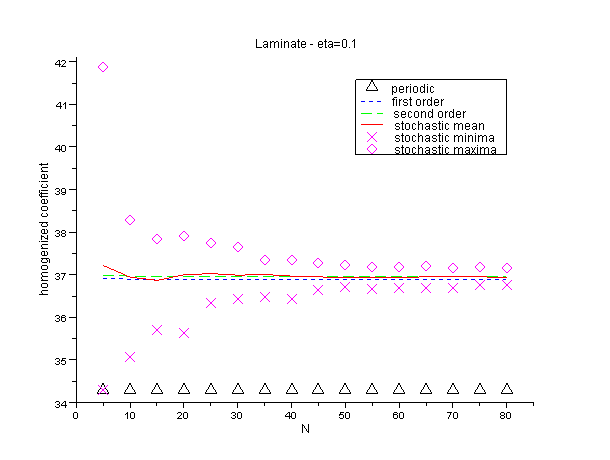} 
\includegraphics[width=14cm, height=10cm]{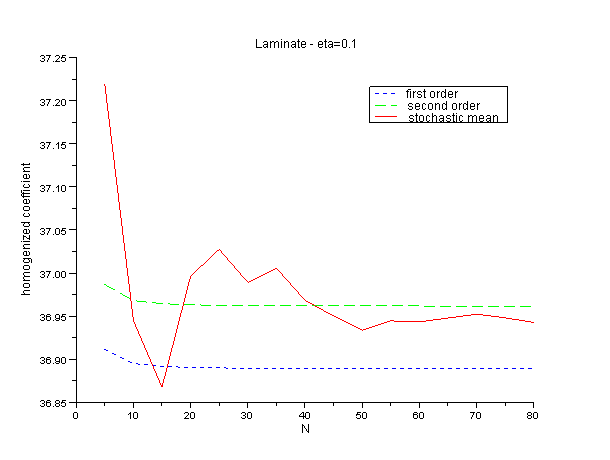} 
\caption{Results for material 2 and $\eta=0.1$. Above: complete results. Below: zoom on $A_{\eta}^{*,N}$ and the first and second-order corrections.}
\label{convlamine01}
\end{figure}

\begin{figure}[H]
\center
\includegraphics[width=14cm, height=11cm]{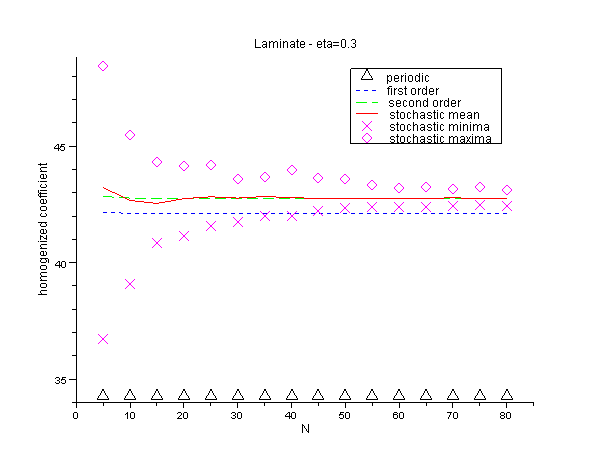} 
\includegraphics[width=14cm, height=11cm]{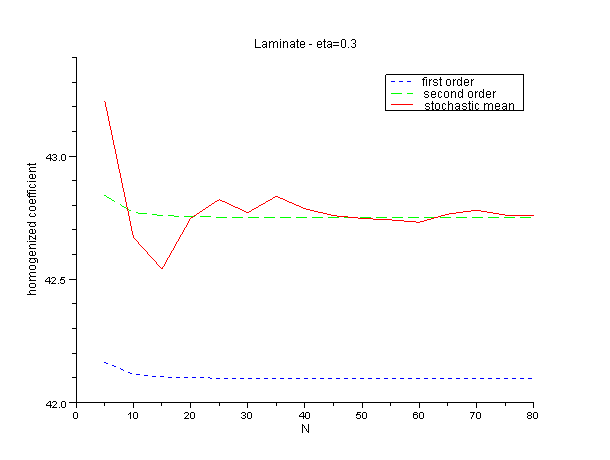} 
\caption{Results for material 2 and $\eta=0.3$. Above: complete results. Below: zoom on $A_{\eta}^{*,N}$ and the first and second-order corrections.}
\label{convlamine03}
\end{figure}

\begin{figure}[H]
\center
\includegraphics[width=14cm, height=11cm]{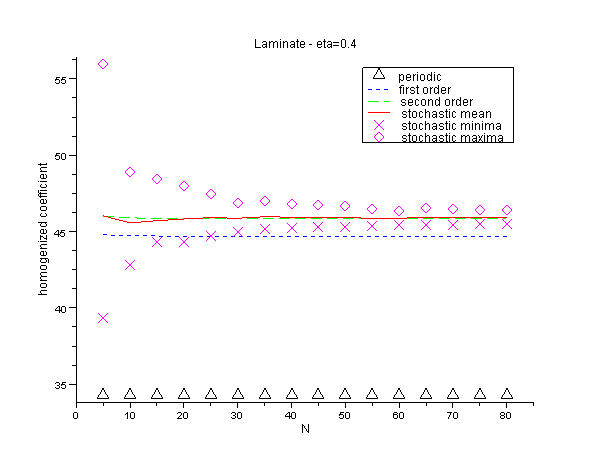} 
\includegraphics[width=14cm, height=11cm]{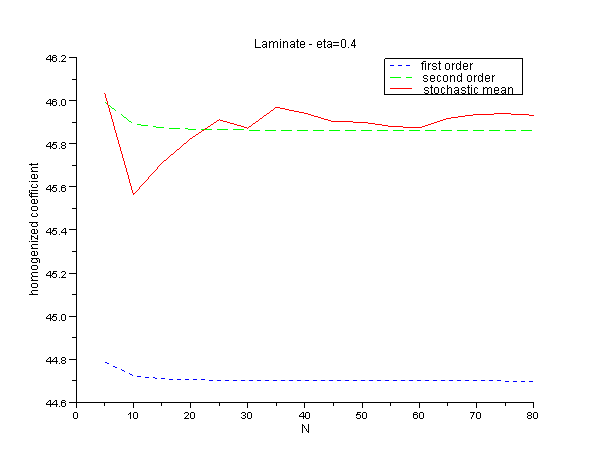} 
\caption{Results for material 2 and $\eta=0.4$. Above: complete results. Below: zoom on $A_{\eta}^{*,N}$ and the first and second-order corrections.}
\label{convlamine04}
\end{figure}

We now show for material $2$ the results for $\eta=0.1$, $\eta=0.3$ and $\eta=0.4$ (Figures \ref{convlamine01}, \ref{convlamine03} and \ref{convlamine04} respectively).

\newpage

\newpage

\newpage




\subsubsection{Comments} \label{comments}

Notice on the results for both materials (it is especially clear on the close-ups) that the first and second-order corrections $A_1^{*,N}$ and $A_2^{*,N}$ converge very fast in function of $N$, and in particular, as expected, much faster than the stochastic computation. Convergence of these deterministic computations is actually typically reached for $N=10$.\\

Then, for all values of $\eta$, it is clear that the first-order correction enables to get substantially closer to $A_{\eta}^*$. The interest of the second-order term is also obvious as $\eta$ gets larger, and we stress that the results are still excellent for $\eta$ as large as $0.5$, so that our approach is robust.\\

It is interesting to get some insight on the rate of convergence of the first-order correction, and to see whether the theoretical results of Lemma \ref{rate} still hold beyond the somewhat restrictive assumptions set in this lemma ($d\geq 3$, piecewise Hölder regularity on $A_{per}$ and Dirichlet boundary conditions on $\partial I_N$). Recall that $d$ is equal to $2$ in our tests, and that $A_1^{*,N}$ is computed with periodic boundary conditions on the supercell $I_N$. Moreover, while the lattice of inclusions is piecewise Hölder continuous in the sense of Lemma \ref{rate} (meaning that there is an inclusion stricly contained in the unit cell $Q$ and that the matrix $A_{per}$ is Hölder continuous in each phase), the laminate is not.\\

\begin{figure}[H]
\center
\includegraphics[width=13cm, height=9cm]{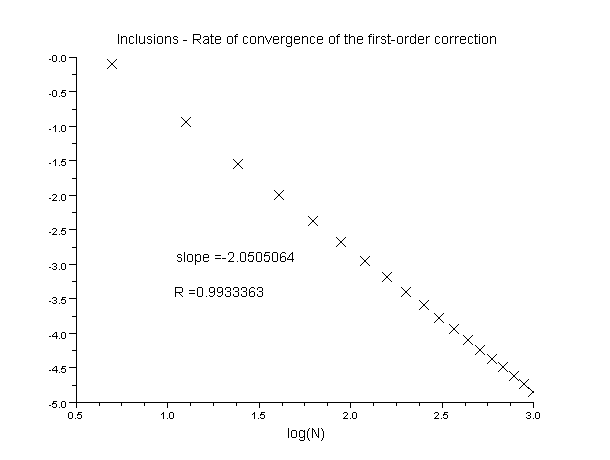} 
\caption{Rate of convergence of the first-order correction for material 1.}
\label{rategraph1}
\end{figure}

\begin{figure}[H]
\center
\includegraphics[width=13cm, height=9cm]{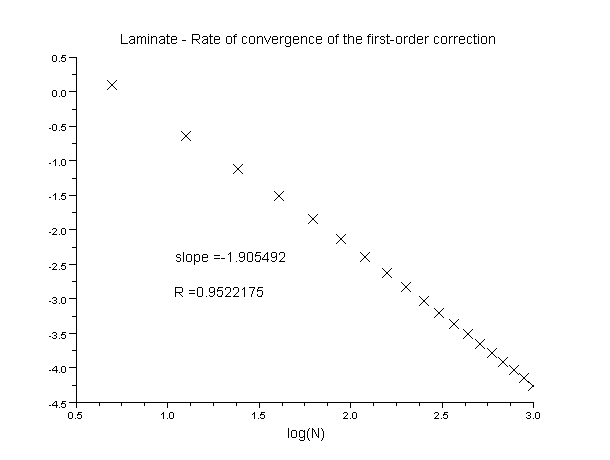} 
\caption{Rate of convergence of the first-order correction for material 2.}
\label{rategraph2}
\end{figure}

 We thus plot, for $N$ going from $1$ to $20$ and for both materials, $\log(|(A_1^{*,N} -A_1^{*})e_1\cdot e_1|)$ in function of $\log(N)$. We recall that $A_1^{*}$ is numerically given by $A_1^{*,80}$. For both materials the $20$ points are arranged in a straight line (Figures \ref{rategraph1} and \ref{rategraph2}). This leads us to perform a linear regression in order to obtain the slope of the lines. As regards material $1$, we find a slope of $-2.05$  and a coefficient of correlation $R = 0.99$. For material $2$, the slope is $-1.9$ with a coefficient of correlation equal to $0.95$. The rate of convergence for both materials is then approximately $\mathcal{O}(N^{-d})$ with $d=2$, which seems to indicate that the result of Lemma \ref{rate} still holds true in these circumstances.


\section{Appendix}

The purpose of this appendix is two-fold. In Section \ref{1dc} we prove that the approach exposed in Section \ref{model}, which relies on formal considerations for general dimensions, is rigorous in dimension $d=1$. In Section \ref{technical}, we prove for convenience of the reader some technical results used in Section \ref{model}. 

\subsection{One-dimensional computations} \label{1dc}

Although we are aware that homogenization theory is very specific in dimension $1$, and can be somehow misleading by its simplicity, it is still important to check that our approach is rigorously founded in this setting. This is the aim of this section.\\

To stress that we work in dimension one, we use lower-case letters $a_{per}$ and $c_{per}$ instead of $A_{per}$ and $C_{per}$, respectively, as well as for all the tensors manipulated.\\

We recall that in dimension one, $a_{per}^{*}$ and $a_{\eta}^*$ are given by the explicit expressions
$$a_{per}^* = \left(\int_{-\frac{1}{2}}^{\frac{1}{2}} \frac{1}{a_{per}}\right)^{-1}, \quad   a_{\eta}^{*} = \left(\mathbb{E}\int_{-\frac{1}{2}}^{\frac{1}{2}} \frac{1}{a_{per} + b_{\eta} c_{per}}\right)^{-1}.$$
This enables us to prove the following elementary result which shows that our approach is correct in dimension one:
\begin{Prop} \label{preuve1d}
In dimension $d=1$, it holds
$$a_{\eta}^{*} = a_{per}^* + \eta \bar{a}_1^{*} + \eta^2 \bar{a}_2^{*} + \mathcal{O}(\eta^3),$$
where $\bar{a}_1^{*}$ and $\bar{a}_2^{*}$ are the limits as $N \rightarrow \infty$ of $a_1^{*,N}$ and $a_2^{*,N}$ defined generally by (\ref{terma1}) and (\ref{terma2}) respectively.
\end{Prop} 

\begin{proof}
We compute
\begin{eqnarray*}
(a_{\eta}^{*})^{-1} &=& \int_{-\frac{1}{2}}^{\frac{1}{2}} \left(\frac{1-\eta}{a_{per}} + \frac{\eta}{a_{per}+c_{per}}\right) \\
&=& \int_{-\frac{1}{2}}^{\frac{1}{2}} \frac{1}{a_{per}} + \eta \int_{-\frac{1}{2}}^{\frac{1}{2}} \left(\frac{1}{a_{per}+c_{per}} - \frac{1}{a_{per}}\right) \\
&=& (a_{per}^{*})^{-1} \left(1-\eta a_{per}^{*} \int_{-\frac{1}{2}}^{\frac{1}{2}} \frac{c_{per}}{a_{per} (a_{per}+c_{per})} \right).
\end{eqnarray*}

This yields the expansion
\begin{eqnarray*} 
a_{\eta}^{*} = a_{per}^{*} + \eta (a_{per}^{*})^2\int_{-\frac{1}{2}}^{\frac{1}{2}} \frac{c_{per}}{a_{per} (a_{per}+c_{per})}+\eta^2 (a_{per}^{*})^3 \left(\int_{-\frac{1}{2}}^{\frac{1}{2}} \frac{c_{per}}{a_{per} (a_{per}+c_{per})}\right)^2\\ 
+ \eta^3 (a_{per}^{*})^4 \left(\int_{-\frac{1}{2}}^{\frac{1}{2}} \frac{c_{per}}{a_{per} (a_{per}+c_{per})}\right)^3 \left(1-\eta a_{per}^{*} \int_{-\frac{1}{2}}^{\frac{1}{2}} \frac{c_{per}}{a_{per} (a_{per}+c_{per})} \right)^{-1}\\
= a_{per}^{*} + \eta (a_{per}^{*})^2\int_{-\frac{1}{2}}^{\frac{1}{2}} \frac{c_{per}}{a_{per} (a_{per}+c_{per})}+\eta^2 (a_{per}^{*})^3 \left(\int_{-\frac{1}{2}}^{\frac{1}{2}} \frac{c_{per}}{a_{per} (a_{per}+c_{per})}\right)^2\\ 
+ \eta^3 (a_{per}^{*})^3 \left(\int_{-\frac{1}{2}}^{\frac{1}{2}} \frac{c_{per}}{a_{per} (a_{per}+c_{per})}\right)^3 a_{\eta}^{*}.
\end{eqnarray*}

It follows from (\ref{Lebesgue}) and (\ref{expect}) that the function $\eta \rightarrow a_{\eta}^{*}$ is bounded on $[0,1]$. Therefore
\begin{eqnarray} \label{devexact}
\begin{aligned}
a_{\eta}^{*} =& a_{per}^{*} + \eta (a_{per}^{*})^2\int_{-\frac{1}{2}}^{\frac{1}{2}} \frac{c_{per}}{a_{per} (a_{per}+c_{per})}\\
& +\eta^2 (a_{per}^{*})^3 \left(\int_{-\frac{1}{2}}^{\frac{1}{2}} \frac{c_{per}}{a_{per} (a_{per}+c_{per})}\right)^2 + \mathcal{O}(\eta^3).
\end{aligned}
\end{eqnarray}

We now devote the rest of the proof to verifying that the coefficients of $\eta$ and $\eta^2$ in (\ref{devexact}) are indeed obtained as the limit as $N \rightarrow \infty$ of $a_1^{*,N}$ and $a_2^{*,N}$  generally  defined by (\ref{terma1}) and (\ref{terma2}) respectively, in this particular one-dimensional setting.\\

The one-defect supercell solution $w^{1,0,N}$ generally defined  by (\ref{onedef}) satisfies here
\begin{eqnarray*} 
\left  \{
\begin{aligned}
& -\frac{d}{dx} \left( a_{1}^0 \left(\frac{d}{dx} w^{1,0,N} + 1 \right)\right) = 0 \quad \mathrm{in} \; ]-\frac{N}{2},\frac{N}{2}[,\\
& w^{1,0,N}\, \, N-\mathrm{periodic}.
\end{aligned}
\right.
\end{eqnarray*}
We easily compute
\begin{eqnarray*}
\begin{aligned}
a_1^0 (\frac{d}{dx} w^{1,0,N}+1) &= N \left(\int_{-\frac{N}{2}}^{\frac{N}{2}} \frac{1}{a_{per}+\mathds{1}_{[-\frac{1}{2},\frac{1}{2}]} c_{per}}\right)^{-1} \\
&= N \left(N (a_{per}^*)^{-1}-\int_{-\frac{1}{2}}^{\frac{1}{2}} \frac{c_{per}}{a_{per} (a_{per}+c_{per})} \right)^{-1} \\
& \hspace{-2cm }= a_{per}^{*}+ \frac{(a_{per}^{*})^2}{N}\int_{-\frac{1}{2}}^{\frac{1}{2}} \frac{c_{per}}{a_{per} (a_{per}+c_{per})} + \frac{(a_{per}^{*})^3}{N^2 }\int_{-\frac{1}{2}}^{\frac{1}{2}}\left(\frac{c_{per}}{a_{per} (a_{per}+c_{per})}\right)^2 + o(N^{-2}).
\end{aligned}
\end{eqnarray*}

Thus $a_1^{*,N}$ defined generally by (\ref{terma1}) takes here the form 
\begin{eqnarray*} 
a_1^{*,N} = \int_{-\frac{N}{2}}^{\frac{N}{2}} a_1^0 (\frac{d}{dx} w^{1,0,N}+1) - Na_{per}^{*} =  (a_{per}^{*})^2 \int_{-\frac{1}{2}}^{\frac{1}{2}} \frac{c_{per}}{a_{per} (a_{per}+c_{per})} +o(1),
\end{eqnarray*}
and
\begin{eqnarray} \label{1da1}
a_1^{*,N} \underset{N \rightarrow \infty}{\rightarrow} \bar{a}_1^{*} = (a_{per}^{*})^2 \int_{-\frac{1}{2}}^{\frac{1}{2}} \frac{c_{per}}{a_{per} (a_{per}+c_{per})}. 
\end{eqnarray}

Likewise, for $k \in \llbracket -\frac{N-1}{2}, \frac{N-1}{2} \rrbracket$, 
\begin{eqnarray*}
\begin{aligned}
a_2^{0,k} (\frac{d}{dx} w^{2,0,k,N}+1) &= N \left(\int_{-\frac{N}{2}}^{\frac{N}{2}} \frac{1}{a_{per}+\mathds{1}_{[-\frac{1}{2},\frac{1}{2}]\cup [k-\frac{1}{2},k+\frac{1}{2}]}c_{per}}\right)^{-1} \\
&= N \left(N (a_{per}^*)^{-1}-2\int_{-\frac{1}{2}}^{\frac{1}{2}} \frac{c_{per}}{a_{per} (a_{per}+c_{per})} \right)^{-1} \\
& \hspace{-3cm }=a_{per}^{*}+ 2\frac{(a_{per}^{*})^2}{N}\int_{-\frac{1}{2}}^{\frac{1}{2}} \frac{c_{per}}{a_{per} (a_{per}+c_{per})} + 4\frac{(a_{per}^{*})^3}{N^2 }\int_{-\frac{1}{2}}^{\frac{1}{2}}\left(\frac{c_{per}}{a_{per} (a_{per}+c_{per})}\right)^2 + o(N^{-2}),
\end{aligned}
\end{eqnarray*}

which is independent of $k$ (and so of the distance between the two defects). Hence, $a_2^{*,N}$ defined generally by (\ref{terma2}) writes here
\begin{eqnarray*} 
a_2^{*,N} = (a_{per}^{*})^3 \left(\int_{-\frac{1}{2}}^{\frac{1}{2}} \frac{c_{per}}{a_{per} (a_{per}+c_{per})}\right)^2 + o(1),
\end{eqnarray*}
and
\begin{eqnarray} \label{1da2}
a_2^{*,N}  \underset{N \rightarrow \infty}{\rightarrow} \bar{a}_2^{*} = (a_{per}^{*})^3 \left(\int_{-\frac{1}{2}}^{\frac{1}{2}} \frac{c_{per}}{a_{per} (a_{per}+c_{per})}\right)^2.
\end{eqnarray}

Using (\ref{devexact}), (\ref{1da1}) and (\ref{1da2}), we verify that
$$ a_{\eta}^{*} = a_{per}^{*} + \eta \bar{a}_1^{*}  + \eta^2 \bar{a}_2^{*} + \mathcal{O}(\eta^3).$$ 
\end{proof}

\begin{remark}
The fact that the distance between two defects does not play a role in the computation of $a_2^{*,N}$ is of course specific to the one-dimensional setting. As we have seen, this is not true in higher dimensions where the geometry comes into play.
\end{remark}

\subsection{Some technical lemmas} \label{technical}

The second part of this appendix is different in nature. We prove here three technical lemmas that are useful for our proofs in Section \ref{model}. These results, or related ones, are probably well known and part of the mathematical literature. We prove them here under specific assumptions for the convenience of the reader and for consistency. We acknowledge several instructive discussions with Xavier Blanc on the content of this section.\\

We recall that $Q=[-\frac{1}{2},\frac{1}{2}]^d$ and $I_N = [-\frac{N}{2},\frac{N}{2}]^d$.\\ 

\begin{lemma} \label{conv1}
Consider $f \in L^2(Q)$, and a tensor field $A$ from $\RR^d$ to $\RR^{d\times d}$ such that there exist $\lambda>0$ and $\Lambda > 0$ such that
$$ \forall \xi \in \RR^d, \mathrm{\; a.e \; in \;} x \in \RR^d, \; \lambda |\xi|^2 \leq A(x)\xi\cdot \xi \; \mathrm{and} \; |A(x) \xi| \leq \Lambda |\xi|.$$ 
 Consider $q^N$ solution to 
\begin{eqnarray} \label{perturblem}
\left  \{
\begin{aligned}
& -\mathrm{div} \left(A \nabla q^{N} \right) = \mathrm{div} (\mathds{1}_Q f) \quad \mathrm{in} \; I_N,\\
& q^{N}\,(N \ZZ)^d-\mathrm{periodic}.
\end{aligned}
\right.
\end{eqnarray} 
Then $ \mathds{1}_{I_N} \nabla q^N$ converges in $L^2(\RR^d)$, when $N$ goes to infinity, to $\nabla q^{\infty}$, where $q^{\infty}$ is a~$L^2_{loc}(\RR^d)$ function solving
\begin{eqnarray} \label{perturblemdef}  
\left  \{
\begin{aligned}
& -\mathrm{div} \left(A \nabla q^{\infty} \right) = \mathrm{div} (\mathds{1}_Q f) \quad \mathrm{in} \; \RR^d,\\
& \nabla q^{\infty} \in L^2(\RR^d).
\end{aligned}
\right.
\end{eqnarray} 
\end{lemma}

\begin{proof}

We first obtain a bound on $\|\nabla q^{N} \|_{L^2(I_N)}$ and then, by compactness, extract a limit of this sequence.\\

Multiplying the first line of (\ref{perturb}) by $q^{N}$ and integrating by parts yields
\begin{eqnarray} \label{varia1}
\int_{I_N} A \nabla q^N \cdot \nabla q^N = - \int_Q f \cdot \nabla q^N, 
\end{eqnarray}
from which we deduce
\begin{eqnarray} \label{bound1}
\|\nabla q^{N} \|_{L^2(I_N)} \leq \frac{1}{\lambda}\|f\|_{L^2(Q)}.
\end{eqnarray}

Consider now a bounded domain $\mathcal{D} \subset \mathbb{R}^d$. For $N$ sufficiently large, we have $\mathcal{D} \subset I_N$ and so
$$ \|\nabla q^{N} \|_{L^2(\mathcal{D})} \leq \frac{1}{\lambda}\|f\|_{L^2(Q)}.$$

Thus $\nabla q^{N}$ is bounded in $L^2(\mathcal{D})$ for every bounded subset $\mathcal{D} \subset \mathbb{R}^d$.\\
 
Using diagonal extraction and the weak compactness of $L^2_{loc}(\RR^{d})$, we can classically find a subsequence of $\nabla q^{N}$ such that, without changing the notation for simplicity,
\begin{eqnarray} \label{convgrad}
\begin{aligned}
\nabla q^{N} \rightharpoonup h \quad \mathrm{weakly \; in\; } L^2_{loc}(\RR^{d}).
\end{aligned} 
\end{eqnarray}
We deduce from (\ref{bound1}) and (\ref{convgrad}) that for every bounded subset $\mathcal{D} \subset \mathbb{R}^d$, 
$$\|h\|_{L^2(\mathcal{D} )} \leq \frac{1}{\lambda}\|f\|_{L^2(Q)}.$$
This implies that the vector $h$ is in $L^2(\RR^{d})$.\\

We also deduce from (\ref{convgrad}) that for all $(i,j) \in \llbracket 1, d \rrbracket^2$, $\frac{\partial h_j}{\partial x_i} =  \frac{\partial h_i}{\partial x_j}$. This implies that $h$ is the gradient of a function we call $q^{\infty}$. Since $h \in L^2(\RR^{d})$, $\nabla q^{\infty}=h$ is in $L^2(\RR^{d})$ and $q^{\infty}$ in $L^2_{loc}(\RR^{d})$.\\

Finally, (\ref{convgrad}) yields that $\nabla q^{N}$ converges to $\nabla q^{\infty}$ in $\mathcal{D}'(\RR^d)$. We can then pass to the limit $N \rightarrow \infty$ in the first line of (\ref{perturblem}) and obtain
$$-\mathrm{div} \left(A \nabla q^{\infty} \right) = \mathrm{div} (\mathds{1}_Q f) \quad \mathrm{in} \; \RR^d$$
in the sense of distributions.\\

We have proved that $\nabla q^{N}$ converges up to extraction and weakly in $L^2_{loc}(\RR^d)$ to $\nabla q^{\infty}$, where $q^{\infty}$ is in $L^2_{loc}(\RR^d)$ and solves
\begin{eqnarray} \label{perturblemdefaux}  
\left  \{
\begin{aligned}
& -\mathrm{div} \left(A \nabla q^{\infty} \right) = \mathrm{div} (\mathds{1}_Q f) \quad \mathrm{in} \; \RR^d,\\
& \nabla q^{\infty} \in L^2(\RR^d).
\end{aligned}
\right.
\end{eqnarray} 
We deduce from Lemma \ref{laxlem} thereafter that (\ref{perturblemdefaux}) has a solution unique up to an additive constant, so that $\nabla q^{\infty}$ is uniquely defined. A classical compactness argument then yields that the whole sequence $\nabla q^{N}$ converges weakly to $\nabla q^{\infty}$ in $L^2_{loc}(\RR^d)$.\\

It is clear from what precedes that 
\begin{eqnarray} \label{convfaible}
\begin{aligned}
\mathds{1}_{I_N} \nabla q^N \rightharpoonup \nabla q^{\infty} \quad \mathrm{weakly \; in\; } L^2(\RR^{d}).
\end{aligned} 
\end{eqnarray}

We now prove that the sequence $\mathds{1}_{I_N} \nabla q^N$ actually converges strongly to $\nabla q^{\infty}$ in $L^2(\RR^d)$.\\

Using a cut-off technique as in the proof of Lemma \ref{laxlem} thereafter, we deduce from (\ref{perturblemdef}) that
\begin{eqnarray} \label{varia2}
\int_{\RR^d} A \nabla q^{\infty} \cdot \nabla q^\infty = - \int_Q f \cdot \nabla q^\infty. 
\end{eqnarray}

The weak convergence of $\nabla q^N$ to $\nabla q^{\infty}$ implies that the right-hand side of (\ref{varia1}) converges to the right-hand side of (\ref{varia2}). Consequently,
\begin{eqnarray}\label{conva}
 \int_{I_N} A \nabla q^N \cdot \nabla q^N \rightarrow \int_{\RR^d} A \nabla q^{\infty} \cdot \nabla q^\infty,
\end{eqnarray}
and, denoting by $A_s$ the symmetric part of $A$, (\ref{conva}) is equivalent to
\begin{eqnarray} \label{convas}
 \int_{I_N} A_s \nabla q^N \cdot \nabla q^N \rightarrow \int_{\RR^d} A_s \nabla q^{\infty} \cdot \nabla q^\infty.
\end{eqnarray}

$A_s$ is of course a uniformly coercive tensor field, we can thus define its square root $A_s^{1/2}$. It follows from (\ref{convas}) that
\begin{eqnarray} \label{convnorm}
\begin{aligned}
\|A_s^{1/2} \mathds{1}_{I_N} \nabla  q^N\|_{L^2(\RR^d)} \rightarrow \|A_s^{1/2} \nabla q^{\infty}\|_{L^2(\RR^d)}.
\end{aligned} 
\end{eqnarray}

On the other hand, multiplying (\ref{convfaible}) by $A_s^{1/2}$, we obtain
\begin{eqnarray} \label{convfaiblea}
\begin{aligned}
A_s^{1/2} \mathds{1}_{I_N} \nabla q^N \rightharpoonup A_s^{1/2} \nabla q^{\infty} \quad \mathrm{weakly \; in\; } L^2(\RR^{d}).
\end{aligned} 
\end{eqnarray}

Because of the uniform convexity of $L^2(\RR^{d})$, it is well known that (\ref{convnorm}) and (\ref{convfaiblea}) imply 
\begin{eqnarray} \label{convfortea}
\begin{aligned}
A_s^{1/2} \mathds{1}_{I_N} \nabla q^N \rightarrow A_s^{1/2} \nabla q^{\infty} \quad \mathrm{strongly \;in\;} L^2(\RR^{d}).
\end{aligned} 
\end{eqnarray}

Multiplying (\ref{convfortea}) by $A_s^{-1/2}$, we finally have
\begin{eqnarray} \label{convforte}
\begin{aligned}
\mathds{1}_{I_N} \nabla q^N \rightarrow  \nabla q^{\infty} \quad \mathrm{strongly \;in\;} L^2(\RR^{d}).
\end{aligned} 
\end{eqnarray}

\end{proof}

\begin{lemma} \label{laxlem}
Let $A$ be a tensor field from $\RR^d$ to $\RR^{d\times d}$ such that there exist $\lambda>0$ and $\Lambda > 0 $ such that
$$ \forall \xi \in \RR^d, \mathrm{\; a.e \; in \;} x \in \RR^d, \; \lambda |\xi|^2 \leq A(x)\xi\cdot \xi \; \mathrm{and} \; |A(x) \xi| \leq \Lambda |\xi.$$
Consider $u \in L^2_{loc}(\RR^d)$ solving
\begin{eqnarray} \label{lax} 
\left  \{
\begin{aligned}
& -\mathrm{div} \left(A\nabla u \right) = 0 \quad \mathrm{in} \; \RR^d,\\
& \nabla u\in L^2(\RR^d). 
\end{aligned}
\right.
\end{eqnarray} 
Then $u$ is constant.
\end{lemma}

\begin{proof}
We define a smooth cut-off function $\chi \in C^{\infty}(\RR^d)$ such that $\chi = 1$ in the ball $B_R$, $\chi = 0$ in $\RR^d \backslash B_{2R}$ and $\|\nabla \chi \|_{L^{\infty}(\RR^d)} \leq 2/R$.\\

Multiplying the first line of (\ref{lax}) by $\chi u$ and integrating by parts, we obtain
$$ \int_{\RR^d} A \nabla u \cdot (\nabla u) \, \chi = - \int_{\RR^d} A \nabla u \cdot (\nabla \chi) \, u.$$
Using the Cauchy-Schwarz inequality, this yields
\begin{eqnarray} \label{firstineq}
\begin{aligned}
\int_{B_R} |\nabla u|^2 &\leq  \frac{\Lambda}{\lambda} \|\nabla \chi \|_{L^{\infty}(\RR^d)} \left(\int_{B_{2R}\backslash B_R} |\nabla u|^2\right)^{1/2} \left(\int_{B_{2R}\backslash B_R} |u|^2\right)^{1/2}  \\
&\leq  \frac{2\Lambda}{R \lambda} \left(\int_{B_{2R}\backslash B_R} |\nabla u|^2\right)^{1/2} \left(\int_{B_{2R}\backslash B_R} |u|^2\right)^{1/2}.
\end{aligned} 
\end{eqnarray}

Defining
$$u_R = \frac{1}{|B_{2R}\backslash B_R|} \int_{B_{2R}\backslash B_R} u,$$
it is clear that $u-u_R$ is also a solution to (\ref{lax}) so that the above computations are valid for $u-u_R$. Since $\nabla (u-u_R) = \nabla u$, we deduce from (\ref{firstineq}) that
\begin{eqnarray} \label{secondineq}
\int_{B_R} |\nabla u|^2 \leq  \frac{2\Lambda}{R \lambda} \left(\int_{B_{2R}\backslash B_R} |\nabla u|^2\right)^{1/2} \left(\int_{B_{2R}\backslash B_R} |u-u_R|^2\right)^{1/2}.  
\end{eqnarray}

We next apply the Poincaré-Wirtinger inequality to $u-u_R$ on $B_{2R}\backslash B_R$. There exists a constant $C(R)$ which depends only on $R$ such that
$$  \int_{B_{2R}\backslash B_R} |u-u_R|^2 \leq C(R) \int_{B_{2R}\backslash B_R} |\nabla u|^2.$$

An easy scaling argument shows that $C(R)$ is equal to $R$ times the Poincaré-Wirtinger constant on $B_2 \backslash B_1$, so that there exists a constant $C$ such that 
\begin{eqnarray} \label{PWball}
\int_{B_{2R}\backslash B_R} |u-u_R|^2 \leq C R \int_{B_{2R}\backslash B_R} |\nabla u|^2.
\end{eqnarray}

We deduce from (\ref{secondineq}) and (\ref{PWball}) that
\begin{eqnarray} \label{absurde}
\int_{B_R} |\nabla u|^2 \leq  \frac{2 C \Lambda}{\lambda} \int_{B_{2R}\backslash B_R} |\nabla u|^2.  
\end{eqnarray}

Since $\nabla u \in L^2(\RR^d)$, the left-hand side of (\ref{absurde}) converges to $\int_{\RR^d} |\nabla u|^2$ when $R \rightarrow \infty$, and the right-hand side of (\ref{absurde}) converges to $0$. Then $\nabla u = 0$ and $u$ is a constant.
\end{proof}

\begin{lemma} \label{regdec}
For $d \geq 3$, consider a $\ZZ^d$-periodic tensor field $A$ such that there exist $\lambda>0$ and $\Lambda > 0 $ such that
$$ \forall \xi \in \RR^d, \mathrm{\; a.e \; in \;} x \in \RR^d, \; \lambda |\xi|^2 \leq A(x)\xi\cdot \xi \; \mathrm{and} \; |A(x) \xi| \leq \Lambda |\xi.$$
Assume that the unit cell $Q$ contains an inclusion $D$, the boundary of which has regularity $\mathcal{C}^{1,\mu}$ for some $0 < \mu < 1$, and such that $\mathrm{dist}(D,\partial Q)>0$. Assume also that $A_{per}$ is Hölder continuous in $\overline{D}$ and in $\overline{Q \backslash D}$.\\

Let  $f$ be a function in $L^2(Q)$.\\

There exists a unique solution $u \in L^2_{loc}(\RR^d)$ to
\begin{eqnarray} \label{unbounded}
\left  \{
\begin{aligned}
& -\mathrm{div} \left(A \nabla u \right) = \mathrm{div} \left(\mathds{1}_Q f \right) \quad \mathrm{in} \; \RR^d,\\
& \nabla u\in L^2(\RR^d), \lim_{|x| \rightarrow \infty} u(x) = 0.
\end{aligned}
\right.
\end{eqnarray}  

Defining also $u_0$ the unique solution to
\begin{eqnarray} \label{bounded} 
\left  \{
\begin{aligned}
& -\mathrm{div} \left(A \nabla u_0 \right) = \mathrm{div}\left(\mathds{1}_Q f \right) \quad \mathrm{in} \; \mathcal{O}, \\
& u_0 \in H^1_0(\mathcal{O}),
\end{aligned}
\right.
\end{eqnarray}
where $\mathcal{O}$ is a bounded domain of $\RR^d$ containing $Q$ and such that $\mathrm{dist}(\partial \mathcal{O}, Q) >1$, there exists a constant $K$ which depends only on $\lambda$, $\Lambda$, $\mu$, $d$, $f$ and the Hölder exponents, and \emph{not} on the domain, such that for $|x| \geq 1$, it holds
$$|u_0(x)| \leq \frac{K}{|x|^{d-1}}, \quad |\nabla u_0(x)| \leq \frac{K}{|x|^{d}},$$
$$|u(x)| \leq \frac{K}{|x|^{d-1}}, \quad |\nabla u(x)| \leq \frac{K}{|x|^{d}}.$$

\end{lemma}

\begin{proof}
Let $G_0$ be the Green kernel associated with $A$ with homogeneous Dirichlet boundary conditions on $\mathcal{\partial O}$, uniquely defined by
\begin{eqnarray*}
\left  \{
\begin{aligned}
& - \mathrm{div}(A \nabla G_0(\cdot,y)) = \delta_y \quad \mathrm{in} \; \mathcal{O},\\
& G_0(\cdot,y) \in W_0^{1,1}(\mathcal{O}),
\end{aligned}
\right.
\end{eqnarray*}
and $G$ be the Green kernel associated with $A$ on $\RR^d$, unique solution to
\begin{eqnarray*}
\left  \{
\begin{aligned}
& - \mathrm{div}(A \nabla G(\cdot,y)) = \delta_y \quad \mathrm{in} \; \RR^d,\\
& G(\cdot,y) \in W^{1,1}_{loc}(\mathbb{R}^d) \cap H^1(\RR^d \backslash B(y,1)).
\end{aligned}
\right.
\end{eqnarray*}
 We deduce from arguments stated in \cite[Lemma 4.2]{Bonnetier} and relying on \cite[Theorem 3.3]{GW}, and on \cite[Lemma 16]{Lin} when $A$ is Hölder continuous and \cite[Theorem 1.9]{Li} when $A$ is piecewise Hölder continuous, that there exists a constant $K$ depending only on $\lambda$, $\Lambda$, $\mu$, $d$ and the Hölder exponents, and not on the domain, such that
\begin{eqnarray} 
\forall (x,y) \in \mathcal{O}, \; \;  |\nabla_y G_0(x,y)| \leq \frac{K}{|x-y|^{d-1}}, \quad |\nabla_x \nabla_y G_0(x,y)| \leq \frac{K}{|x-y|^{d}}, \label{decay0}\\
\forall (x,y) \in \RR^d, \; \;  |\nabla_y G(x,y)| \leq \frac{K}{|x-y|^{d-1}}, \quad |\nabla_x \nabla_y G(x,y)| \leq \frac{K}{|x-y|^{d}}. \label{decay}
\end{eqnarray} 
It is well known that $u_0$ solution to (\ref{bounded}) can be represented as
\begin{eqnarray} \label{u0}
u_0(x) = \int_{\mathcal{O}} G_0(x,y)\mathrm{div}(\mathds{1}_Q f)(y) dy.
\end{eqnarray}

It is also clear that the function $\tilde{u}$ defined by
\begin{eqnarray} \label{tildeu}
 \tilde{u}(x) = \int_{\RR^d} G(x,y)\mathrm{div}(\mathds{1}_Q f)(y) dy
\end{eqnarray}
is a $H^1_{loc}(\RR^d)$ function which satisfies 
$$-\mathrm{div} \left(A\nabla \tilde{u} \right) = \mathrm{div} \left(\mathds{1}_Q f \right)$$
in the sense of distributions.

Integrating by parts in (\ref{u0}) and (\ref{tildeu}) for $x \notin Q$, we have
\begin{eqnarray} \label{gradient}
u_0(x)= \int_{Q} \nabla_y G_0(x,y) \cdot f(y) dy , \; \; \tilde{u}(x)= \int_{Q} \nabla_y G(x,y) \cdot f(y) dy,
\end{eqnarray}
and then
\begin{eqnarray} \label{hessienne}
\nabla u_0(x) = \int_{Q} \nabla_x \nabla_y G_0(x,y) \cdot f(y) dy, \; \; \nabla \tilde{u}(x) = \int_{Q} \nabla_x \nabla_y G(x,y) \cdot f(y) dy.
\end{eqnarray}

Using estimates (\ref{decay0}) and (\ref{decay}) in (\ref{gradient}) and (\ref{hessienne}) respectively, we find that there exists a constant $K$ depending only on $\lambda$, $\Lambda$, $\mu$, $d$, $f$ and the Hölder exponents, and not on the domain, such that for $|x| \geq 1$, we have
\begin{eqnarray}
|u_0(x)| \leq \frac{K}{|x|^{d-1}} \quad \mathrm{and} \quad |\nabla u_0(x)| \leq \frac{K}{|x|^{d}}, \label{decu0} \\
|\tilde{u}(x)| \leq \frac{K}{|x|^{d-1}} \quad \mathrm{and} \quad |\nabla \tilde{u}(x)| \leq \frac{K}{|x|^{d}} \label{dectildeu}.
\end{eqnarray}

The function $\tilde{u}$ being in $H^1_{loc}(\RR^d)$, we deduce from (\ref{dectildeu}) that $\nabla \tilde{u} \in L^2(\RR^d)$. Consequently, $\tilde{u}$ solves
\begin{eqnarray} \label{systildeu}
\left  \{
\begin{aligned}
& -\mathrm{div} \left(A\nabla \tilde{u} \right) = \mathrm{div} \left(\mathds{1}_Q f \right) \quad \mathrm{in} \; \RR^d,\\
& \nabla \tilde{u} \in L^2(\RR^d).
\end{aligned}
\right.
\end{eqnarray}  

We know from Lemma \ref{laxlem} that (\ref{systildeu}) has a solution unique up to an additive constant. It follows from (\ref{dectildeu}) that $\tilde{u}$ converges to zero at infinity, so that $\tilde{u}=u$ unique solution to (\ref{unbounded}).

\end{proof}


\begin{thebibliography}{00}
%

\bibitem{these_Arnaud} A. Anantharaman, {\em Thèse de l'Université Paris-Est}, in preparation.

\bibitem{ALB2} A. Anantharaman, C. Le Bris, {\em Elements of mathematical foundations for a numerical approach for weakly random homogenization problems}, preprint available on this archive.

\bibitem{ALBcras} A. Anantharaman, C. Le Bris, {\em Homogenization of a weakly randomly perturbed periodic material}, C. R. Acad. Sci. Paris Série I, Vol. 348 (9-10) (2010), pp. 529-534.

\bibitem{AG} G. Allaire, S. Gutierrez, {\em Optimal design in small amplitude homogenization}, ESAIM: M2AN, Vol. 41 no.3 (2007), pp. 543-574.

\bibitem{Lin} M. Avellaneda, F.-H. Lin, {\em Compactness methods in the theory of homogenization}, Communications on Pure and Applied Mathematics, Vol. XL (1987), pp. 803-847.

Comput. Methods Appl. Mech. Engrg., 172 (1999), pp. 27-77.

\bibitem{Bonnetier} M. F. Ben Hassen, E. Bonnetier, {\em An asymptotic formula for the voltage potential in a perturbed $\epsilon$-periodic composite medium containing misplaced inclusions of size $\epsilon$}, Proceedings of the Royal Society of Edinburgh, 136A (2006), pp. 669-700.

\bibitem{BLL} X. Blanc, C. Le Bris, P.-L. Lions, {\em Stochastic homogenization and random lattices}, J. Math. Pures Appl., 88
  (2007), pp. 34-63. 

\bibitem{BP} A. Bourgeat, A. Piatnitski, {\em Approximations of effective coefficients in stochastic homogenization}, Annales de l'Institut Henri Poincaré (B) Probabilités et Statistiques, 40 no. 2 (2004), pp. 153-165.

\bibitem{Costaouec} R. Costaouec, C. Le Bris, F. Legoll, {\em Numerical approximation of a class of problems in stochastic homogenization}, C. R. Acad. Sci. Paris Série I, Vol. 348 (1-2) (2010), pp. 99-103.

\bibitem{GW}  M. Grüter, K-O. Widman, {\em The Green function for uniformly elliptic equations}, manuscripta math., 37 (1982), pp. 303-342. 

\bibitem{JKO} V. V. Jikov, S. M. Kozlov, O. A. Oleinik, {\em Homogenization of Differential Operators and Integral Functionals}, Springer Verlag (1994). 

\bibitem{enumath} C. Le Bris, {\em Some numerical approaches for ``weakly'' random homogenization}, to appear in Proceedings of the ENUMATH 2009 Conference, Springer (2010).

\bibitem{Li} Y. Li, L. Nirenberg, {\em Estimates for elliptic systems from composite material}, Communications on Pure and Applied Mathematics, Vol. LVI (2003), pp. 892-925.

\bibitem{Sakata} S. Sakata, F. Ashida, T. Kojima, M. Zako, {\em Three-dimensional stochastic analysis using a perturbation-based homogenization method for elastic properties of composite material considering microscopic uncertainty}, International Journal of Solids and Structures, 45 (2008), pp. 894-907.

\bibitem{Tartar} L. Tartar, {\em H-measures, a new approach for studying homogenisation, oscillations and concentration effects in partial differential equations}, Proceedings of the Royal Society of Edinburgh, 115A (1990), pp. 193-230. 

\bibitem{Matthieu} M. Thomas, {\em Propriétés thermiques de matériaux composites : caractérisation expérimentale et approche microstructurale}, Thèse de l'Université de Nantes, Laboratoire de Thermocinétique, CNRS-UMR 6607 (2008).

\end{thebibliography}
\end{document}